\documentclass[lfeqn,12pt]{article}
\usepackage{amssymb,bbm,pifont,mathtools,mathabx,mathrsfs,%fdsymbol,
}
\usepackage{minitoc}
 
\usepackage[nottoc]{tocbibind}
 
\usepackage{multirow}
\usepackage[multiple]{footmisc}
\usepackage{xcolor}
\usepackage[colorlinks=true]{hyperref}
\newcommand{\giu}{{\medskip\noindent}}

\newcommand{\Giu}{{\bigskip\noindent}}
\newcommand{\nl}{{\smallskip\noindent}}
\newcommand{\noi}{{\noindent}}

\newtheorem{theorem}{Theorem}
\newtheorem{definition}[theorem]{Definition}
\newtheorem{proposition}[theorem]{Proposition}
\newtheorem{lemma}[theorem]{Lemma}
\newtheorem{remark}[theorem]{Remark}
\newtheorem{sublemma}[theorem]{Sublemma}
\newtheorem{corollary}[theorem]{Corollary}
\newtheorem{assumption}[theorem]{Assumption}
\newtheorem{notationalrem}[theorem]{Notational Remark}
\newtheorem{tools}[subsection]{$\negsp\negsp$}
\newcommand\asm[1]{ \begin{assumption}\label{#1} }
\newcommand\easm{ \end{assumption} }
\newcommand\dfn[1]{ \begin{definition}\label{#1} }
\newcommand\dfntwo[2]{ \begin{definition}[#2]\label{#1} }
\newcommand\edfn{ \end{definition} }
\newcommand\rem[1]{ \begin{remark}\label{#1} \small \rm}
\newcommand\remtwo[2]{ \begin{remark}[#2]\label{#1} \rm}
\newcommand\erem{ \end{remark} }
\newcommand\thm[1]{ \begin{theorem}\label{#1}}
\newcommand\thmtwo[2]{ \begin{theorem}[#2]\label{#1}}
\newcommand\ethm{ \end{theorem} }
\newcommand\pro[1]{ \begin{proposition}\label{#1}}       
\newcommand\protwo[2]{ \begin{proposition}[#2]\label{#1}}
\newcommand\epro{ \end{proposition} }
\newcommand\lem[1]{ \begin{lemma}\label{#1}}
\newcommand\lemtwo[2]{ \begin{lemma}[#2]\label{#1}}
\newcommand\elem{ \end{lemma} }
\newcommand\sublem[1]{ \begin{sublemma}\label{#1}}
\newcommand\sublemtwo[2]{ \begin{sublemma}[#2]\label{#1}}
\newcommand\esublem{ \end{sublemma} }
\newcommand\cor[1]{ \begin{corollary}\label{#1}}
\newcommand\cortwo[2]{ \begin{corollary}[#2]\label{#1}}
\newcommand\ecor{ \end{corollary} }
\newcommand\notrem[1]{ \begin{notationalrem}\label{#1} \sl}
\newcommand\enotrem{ \end{notationalrem} }

\newcommand\average[1]{{ \left\langle #1 \right\rangle}}

\newcommand\equ[1]{{\rm (\ref{#1})}}
\newcommand\beq[1]{ \begin{equation}\label{#1} }
\newcommand{\eeq}{ \end{equation} }
\newcommand{\beqno}{ \[ }
\newcommand{\eeqno}{ \] }
\newcommand\beqa[1]{ \begin{eqnarray} \label{#1}}
\newcommand{\eeqa}{ \end{eqnarray} }
\newcommand{\beqano}{ \begin{eqnarray*} }
\newcommand{\eeqano}{ \end{eqnarray*} }

\newcommand{\proof}{\par\medskip\noindent{\bf Proof\ }}

\newcommand{\dst}{\displaystyle}
\newcommand{\qed}{\hskip.5truecm
            \vrule width 1.7truemm height 3.5truemm depth 0.truemm
            \par\Giu}

\newcommand\ovl[1]{ \overline {#1} }

\newcommand\su[1]{ \frac{1}{ {#1}} }

\newcommand\osc{ {\, \rm osc\, }}

\newcommand{\diag}{{ \, \rm diag \, }}

\newcommand{\io}{{\infty }}

\newcommand{\dpr}{ {\partial}   }

\newcommand\eqby[1]{\stackrel{\equ{#1}}{=}}
\newcommand\leby[1]{\stackrel{\equ{#1}}{\le}}

\newcommand\geby[1]{\stackrel{\equ{#1}}{\ge}}
\newcommand\gtby[1]{\stackrel{\equ{#1}}{>}}
\renewcommand{\Im}{{\rm \, Im\,}}

\newcommand{\negsp}{\hspace{-.04truecm}}
\newcommand\ex{\, e}%{\mathbbmss e}}

\renewcommand{\a }{ {\alpha}   }

\newcommand{\g}{ {\gamma}   }

\newcommand{\D}{ {\Delta}   }
\newcommand{\vae }{ {\varepsilon}   }
\newcommand{\e }{ {\varepsilon}   }

\renewcommand{\th }{ \hat \epsilon   }

\newcommand{\torsion}{\theta}

\renewcommand{\k}{ {\kappa}   }
\renewcommand{\l}{ {\lambda}   }

\newcommand{\m}{ {\epsilon}   }
\newcommand{\n}{ {\nu}   }

\newcommand{\p}{ {\pi}   }

\newcommand{\s}{ {\sigma}   }

\renewcommand{\t}{ {\tau}   }
\newcommand{\f}{ {\varphi}   }

\renewcommand{\o}{ {\omega}   }

\newcommand{\torus}{ {\mathbb{ T}}   }
\renewcommand{\natural}{ {\mathbb{ N}}   }
\newcommand{\real}{ {\mathbb{ R}}   }
\newcommand{\integer}{ {\mathbb{ Z}}   }
\newcommand{\complex}{ {\mathbb { C}}   }

\newcommand{\tn}{ {\torus^d} }

\newcommand{\rn}{ {\real^d}   }

\newcommand{\cn}{ {\complex^d }   }
\newcommand{\zn}{ {\integer^d }   }

\newcommand\ppu{{ (1) }}
\newcommand\ppd{{ (2) }}
\newcommand\ppt{{ (3) }}

\font\teneufm=eufm10
\font\seveneufm=eufm7
\font\fiveeufm=eufm5
\newfam\eufmfam
\textfont\eufmfam=\teneufm
\scriptfont\eufmfam=\seveneufm
\scriptscriptfont\eufmfam=\fiveeufm

\newcommand\appA[1]{\section{#1}\label{app:A}
\renewcommand{\theequation}{A.\arabic{equation}}
           \setcounter{equation}{0}
\renewcommand{\thetheorem}{A.\arabic{theorem}}
           \setcounter{theorem}{0}           
           }
\newcommand\appB[1]{\section{#1}\label{app:B}
\renewcommand{\theequation}{B.\arabic{equation}}
           \setcounter{equation}{0}}
\newcommand\appC[1]{\section{#1}\label{app:C}
\renewcommand{\theequation}{C.\arabic{equation}}
           \setcounter{equation}{0}}

\def\bks{{\backslash}}

\def\uno{{\mathbbm 1}}

\def\id{{\rm id }}

\newcommand{\wt}{\widetilde}

\begin{document}

%%%%%%%%%%%%%%%%%%%%%%%%%%%
%%%%%%%   title page

\date{\small \today}

\title{{\bf  V.I.~Arnold's ``pointwise'' KAM Theorem}\\
}
\author{ 
L. Chierchia \& 
C. E. Koudjinan \\
\vspace{-0.25truecm}
{\footnotesize Dipartimento di Matematica\ ,  Universit\`a  ``Roma Tre"}
\\ 
\vspace{-0.25truecm}
{\footnotesize Largo S. L. Murialdo 1, I-00146 Roma (Italy) }
\\
\vspace{-0.25truecm}
{\footnotesize luigi@mat.uniroma3.it,  
ckoudjinan@mat.uniroma3.it}
\\ 
\vspace{-0.25truecm}
}
\maketitle

\begin{abstract}
\noindent
We review V.I.~Arnold's 1963 celebrated paper \cite{ARV63} {\sl Proof of A.N.~Kolmogorov's theorem on the conservation of conditionally periodic motions with a small variation in the Hamiltonian}, and prove that, optimising Arnold's scheme, one can get ``sharp'' asymptotic quantitative conditions (as $\vae\to 0$, $\vae$ being the strength of the perturbation).  All constants involved are explicitly computed.
\end{abstract}
\tableofcontents

\section{Introduction}
\begin{itemize}
\item[\bf  a.]
``{\sl One of the most remarkable of A.N.~Kolmogorov's mathematical achievements is his work on classical mechanics of 1954}'': this is the beginning of  V.I.~Arnold's celebrated paper {\sl Proof of A.N.~Kolmogorov's theorem on the conservation of conditionally periodic motions with a small variation in the Hamiltonian} \cite{ARV63}, published in 1963,  on the occasion of A.N.~Kolmogorov's 60th birthday. Few lines after, Arnold adds: ``{\sl Its deficiency has been that complete proofs have never been published}. \\
Even though one could argue whether Kolmogorov's proof in \cite{Ko1954} is ``complete'' or not (see, e.g., \cite{CL08}), Arnold's paper is certainly a  milestone of modern dynamical systems, which not only contains a {\sl complete and detailed} proof of Kolmogorov's Theorem, but, also, introduces new original, technical ideas, of enormous impact in finite and infinite dimensional systems (for reviews, see, e.g., \cite{AKN06} or \cite{CP19}).   

\item[\bf b.]  
Kolmogorov's 1954 theorem in classical mechanics \cite{Ko1954} (see, also, \cite{CL08}), deals, as is well known, with the persistence, for small $\e$,  of Lagrangian invariant tori of analytic integrable systems governed by a nearly integrable Hamiltonian 
\beq{Horiginal}
H(y,x)=K(y)+\e P(y,x)\ ,
\eeq  
where $(y,x)\in\real^d\times \torus^d$ are  
standard symplectic action--angle variables. In short, the theorem says that:

\nl
{\sl 
for small $\e$, non--degenerate Diophantine unperturbed Lagrangian tori 
persist}

\giu
Let us recall that 
``Diophantine'' means that the unperturbed torus $\mathcal{T}_{\o,0}\coloneqq\{y_0\}\times\torus^d$, which is invariant for the flow $\phi_K^t$ governed by the integrable Hamiltonian $K$, is such that the frequency
$\o\coloneqq K_y(y_0)$ is Diophantine, i.e., it satisfies, for some $\a,\t>0$,  
\beq{nota1}
|\o\cdot k|\coloneqq \sum_{j=1}^d |\o_j k_j|\ge \frac{\a}{|k|^\t}\,,\qquad\forall \ k\in\integer^d\bks\{0\}\,;
\eeq
``non--degenerate'' means that the Hessian of $K$ at $y_0$ is invertible; finally, ``persists'' means that $\mathcal{T}_{\o,0}$ deforms, for positive small enough $\e$,  into a a Lagrangian\footnote{A Lagrangian manifold is a submanifold of dimension
$d$ on which the restriction of the two form $\sum_{j=1}^d dy_j\wedge dx_j$ vanishes.} torus  $\mathcal{T}_{\o,\vae}$ invariant for $\phi_H^t$.

\nl
The scheme on which Arnold's proof of Kolmogorov's theorem is based, while sharing two basic ideas of Kolmogorov's approach -- namely, the use of a quadratic symplectic iterative method and the idea of keeping fixed the Diophantine frequency of the motion --  is quite different from Kolmogorov's scheme in the following respects.\\
 First, for a fixed frequency, Arnold constructs an embedded, Lagrangian invariant torus obtained as a limit of symplectic transformations on action domains shrinking to a single point; in contrast, Kolmogorov conjugates the given Hamiltonian to a complete normal form admitting a Lagrangian invariant torus with the prescribed frequency. 
 \\
 A key difference between these two approaches is that, Arnold, at each step of the iteration, needs to control only a finite number of small divisors\footnote{To work with a finite number of divisors, Arnold introduces a Fourier cut--off (depending, in view of analyticity, logarithmically on the size of the perturbation), an idea which has been widely followed also in infinite dimensional Hamiltonian perturbation theory.}, which however depend on actions (this being the reason for the shrinking to one point of the action domains), while in the denominators appearing in Kolmogorov's scheme there enters only the prefixed Diophantine frequency, allowing one to control at once all small divisors, and also to work with smaller and smaller domains, which contain a fixed open set, allowing one, in the end, 
 to get  a genuine symplectic transformation. \\
 A clever quantitative revisitation of Kolmogorov's scheme (\cite{V08}) shows that such a scheme leads to optimal asymptotic   estimates (as $\e\to 0$). {\sl We shall show below that this is true also  for Arnold's original ``pointwise'' scheme}.

\item[\bf  c.]  
Kolmogorov's and Arnold's schemes are ```pointwise'' in the sense that they deal with the continuation of a single prefixed unperturbed Lagrangian  torus with Diophantine frequency. This is in contrast with  versions of the KAM theorem\footnote{Striclty speaking, there does not exists a {\sl KAM Theorem} (``KAM'' standing for the initials of A.N.~Kolmogorv, V.I.~Arnold and J.K.~Moser), however, normally, it refers to (variations of) Kolmogorov's theorem. Here, we follow this tradition.}
dealing with the persistence of {\sl sets} of simultaneously persistent invariant tori, see 
\cite{ARV63}, \cite{Neish81}, \cite{JP82}, \cite{CG82}. We point out that, actually, Arnold's original formulation of the KAM theorem in \cite{ARV63} belongs to  this second kind of theorem as it states the existence of a set of simultaneously invariant tori, however, {\sl the proof} is pointwise in nature and its scheme is exactly the scheme we follow closely here. Typically, especially when one is concerned with {\sl lower dimensional} invariant tori, it is not possible to construct a single torus with some pre--assigned property, but, rather, one obtains ``Cantor'' families of persistent tori (compare, e.g.,  \cite{CP19}).

\item[\bf d.]
The smallness condition, i.e., how small  the perturbation has to be in order for the perturbed invariant torus to exist, depends on local analytic properties of $K$ (and on the analytic norm of $P$).
In particular,  the main quantitative ``competition'' is between $\vae$ and the size of the small divisors appearing in the iterative scheme, the size of which may be measured by the ``homogeneous Diophantine constant'' $\a$ (compare Eq.~\equ{nota1}) of the prefixed frequency $\o=K_y(y_0)$.
\\
The most important quantitative relations may be easily understood by looking at explicitly solvable examples, i.e., at integrable systems.
\\
To illustrate this point, let us consider, for example, a  simple pendulum with   gravity $\e$,  
\beq{pendulum}
H(y,x)=\frac12 y^2 +\vae (\cos x-1) \ ,
\eeq
viewed as an $\vae$--perturbation of the non--degenerate Hamiltonian $K(y)\coloneqq\frac12 y^2$, (here, $d=1$).
The energy zero level  $\{H=0\}$ corresponds to the separatrix, i.e., 
$$y=\pm \sqrt{2\e (1-\cos x)}\ ,$$
which shows immediately that in the region $\mathcal{S}:=\{|y|\le 2\sqrt\e\}$ there are no homotopically trivial invariant tori (curves) or, equivalently, no Lagrangian invariant curves, which are graphs over the angle 
variable (``primary tori''). In other words, the region of action space where unperturbed curves 
$\{y_0\}\times \torus$ may be continued into invariant Lagrangian invariant curves, which stay out of the ``singular region'' 
$\mathcal{S}$ are such that:
\beq{nec.cond.0}
|y_0|>2\sqrt\e\ .
\eeq
Now, the resonant relations $|K_y(y_0)\cdot k|$ become, in this one--dimensional example,   simply $|y_0| |k|$ and the Diophantine condition is, therefore,  equivalent to requiring that $\a=|y_0|$  (recall \equ{nota1}), and the necessary condition \equ{nec.cond.0} becomes:
\beq{nec.cond}
\frac\e{\a^2}<\frac14\ .
\eeq
Another fact  that can be easily extracted from this example concerns the oscillations of (primary) invariant tori\footnote{A primary Lagrangian torus is a graph over the angles $\{(y,x)|\ y=U(x)\,, x\in \torus^d\}$ and its oscillation is given by $\sup_{x,x'} |U(x)-U(x')|$.}. \\
For $y_0>0$ the invariant (primary) curves are given by 
$$y_\e(x):= \sqrt{y_0^2 +2\e (1-\cos x)}=y_0+v_\e(x)\ ,\qquad $$
with
$$
v_\e(x):= \frac{2\e(1-\cos x)}{y_0+\sqrt{y_0^2+2\e(1-\cos x)}}\ .
$$
Thus, one has that 
$$\osc(y_\e)=\osc (v_\e)\ge v_\e(\p)-v_\e(0)= \frac{4\e}{y_0+\sqrt{y_0^2+4\e}}=\frac\e{y_0}\ \frac{4}{1+\sqrt{1+ 4\e/y_0^2}}\ ,
$$
which, in view of \equ{nec.cond}, yields the relation
\beq{osc.int}
\osc(v_\e)\ge \frac4{1+\sqrt2}\cdot \frac\e\a 
\eeq
Below,   we shall prove that the enhanced Arnold's scheme 
leads to a smallness condition of the type  (compare \equ{smcondwhL} below)
\beq{condA}
\frac{\e}{\a^2} < c\ ,
\eeq
(for an $\e$ and $\a$ independent constant $c$), which is in agreement with \equ{nec.cond}. 
\\
Furthermore, we shall also show that Arnold's scheme leads to a bound on the oscillations of persistent tori given as graphs $\{y=y_0+v_*(x), x\in\torus^d\}$
of the form (compare \equ{est} below)
\beq{osc.int*}
\osc(v_*)\le C\cdot   \frac\e\a\ ,
\eeq
(for an $\e$ and $\a$ independent constant $C$), which, in view of \equ{osc.int},
is seen to be optimal (as far as the dependence upon 
$\e$ and $\a$ is concerned), showing the ``quantitative sharpness'' of Arnold's scheme, on which the proof presented below is based.

\nl
Condition \equ{condA} is also the fundamental quantitative relation needed to evaluate the {\sl measure} of the Kolmogorov's set, i.e., the union (in a prefixed bounded domain) of all  primary  tori. Indeed, \equ{condA} leads to a bound on the Lebesgue measure of the  complement of the Kolmogorov's set by a constant times $\sqrt\e$ (compare \cite{Neish81}, \cite{JP82}), which again, comparing with the simple pendulum \equ{pendulum} -- that has a region (the area enclosed by the separatrix) of measure $16\sqrt\e$ free of primary tori -- is seen to be asymptotically optimal. It has to be remarked, however, that
obtaining such an estimate  is quite delicate and far from trivial (for a more detailed discussion on this point, see \cite{BC18}, \cite{KoudjThesis}, \cite{CKpre}).

\item[\bf e.] As is well known, Arnold's scheme is an iterative  Newton scheme yielding a sequence of  ``renormalised Hamiltonians'' 
$$
H_j:=K_j + \e^{2^j} P_j
$$
so that $H_0=H$ is the given nearly integrable Hamiltonian \equ{Horiginal} and, for any $j$, 
$K_j$ is integrable (i.e., depends only on the action variable $y$), real--analytic in a $r_j$--ball around a point $y_j$ close to $y_0$ and satisfies:
\beq{fund.hyp}
\partial_y K_j(y_j)=\o:=\partial_y K(y_0)\ ,\qquad
\det \partial^2_y K_j(y_j)\neq 0\ ,
\eeq
which means that at each step the frequency is kept fixed and that the integrable Hamiltonian $K_j$  is non--degenerate.
The sequence of Hamiltonians $H_j$ are  conjugated, i.e., $H_{j+1}=H_j\circ \phi_j$, with
$\phi_j$ symplectic, closer and closer to the identity. The persistent torus 
$\mathcal{T}_{\o,\e}$ is then obtained as the limit
$$
\lim_{j\to+\io} 
\phi_0\circ \cdots\phi_{j-1}(y_j,\torus^n)\ .
$$
The symplectic transformations $\phi_j$'s are obtained by solving  the classical Hamilton--Jacobi equation
so as to remove quadratically the order of the perturbation. In doing this one cannot take into account all small divisors (which are dense)  and therefore Arnold introduces a Fourier cut--off $\k_j$, which allows him to deal with a finite number of small divisors. In view of the exponential decay of Fourier coefficients, $\k_j$ can be taken $\sim \big| \log\big(e^{2^j \|P_j\|}\big)\big|$, which introduces a logarithmic correction\footnote{For full details, see \S~\ref{Astep} below, and in particular ``Step 1: Construction of  Arnold's transformation''.}, that does not affect the convergence of the scheme. All this is well known. 
\\
The problem  is to equip the scheme with  ``optimal''  quantitative estimates, which may lead, at the end, to the above sharp asymptotic bounds. This involves careful choices of various parameters entering the scheme (see \S~\ref{sec:iteration}) and,  in particular, 
{\sl it is crucial to treat the first step in a  different way with respect to the remaining steps}: this technical, but important, aspect is explained  in Remark~\ref{steps} below.

\item[\bf f.] V.I. Arnold pointed out that his proof extended with little changes to the {\sl iso--energetically non--degenerate} case, i.e.,  when the energy is prescribed and the unperturbed Hamiltonian satisfies the  condition\footnote{The matrix in \equ{Acond} is a $(d+1)\times(d+1)$ matrix, where the upper right corner $\partial_y K$ has to be interpreted as a column vector, while the lower left corner is a raw vector and the zero is a scalar. The condition expresses the fact the map $(y,\l)\mapsto (\l \partial_yK, K)$ is locally invertible.}
\beq{Acond}
\det 
\begin{pmatrix}\partial^2_y K & \partial_y K\\
\partial_y K & 0
\end{pmatrix}
\Big|_{y=y_0} \neq 0\,.
\eeq
Indeed, it would not be difficult to adapt our improved Arnold's scheme also to the iso--energetically non--degenerate case, proving the sharpness of the asymptotic smallness conditions also in this case.

\item[\bf g.] 
Finally, we mention that the quantitative estimates provided in this paper could be used to improve the (exponentially long) stability time of ``nearly--invariant tori'', introduced in   
\cite{DG96}.

\end{itemize}

\section{Notation and quantitative statement of Arnold's Theorem}

\begin{itemize}

\item[\tiny $\bullet$] For $d\in\natural \coloneqq \{1,2,3,...\}$ and  $x,y\in \complex^d$, we let
$x\cdot y\coloneqq x_1 \bar y_1 +\cdots+x_d \bar y_d$ be the standard inner product;
$|x|_1\coloneqq\dst\sum_{j=1}^d |x_j| $ be the $1$--norm,
and  $|x|\coloneqq \dst\max_{1\le j\le n}|x_j|$ be the sup--norm.

\item[\tiny $\bullet$] $\tn\coloneqq \rn/2\p\zn$ is the  standard $d$--dimensional (flat) torus.

\item[\tiny $\bullet$] $\pi_1\colon \cn\times\cn\ni (y,x)\longmapsto y$ and $\pi_2\colon \cn\times\cn\ni (y,x)\longmapsto x$ are the projections on the first and second component respectively.

\item[\tiny $\bullet$] For $\a>0$, $\t\ge d-1\ge 1$, 
\beq{dio}\D_\a^\t\coloneqq \left\{\o\in \rn:|\o\cdot k|\geq \frac{\a}{|k|_1^\t}, \ \ \forall\ 0\not=k\in \zn\right\},
\eeq
is the set of $(\a,\t)$--Diophantine numbers in $\real^d$. 

\item[\tiny $\bullet$] For $r,s>0$, $y_0\in\complex^d$, we denote:
\beqano
\dst\torus^d_s &\coloneqq  &\left\{x\in \cn: |\Im x|<s\right\}/2\p\zn\,,\\
B_r(y_0)&\coloneqq  &\left\{y\in \real^d: |y-y_0|<r\right\}\,,\qquad (y_0\in\real^d)\,,\\
D_r(y_0)&\coloneqq  &\left\{y\in \cn: |y-y_0|<r\right\}\,,\quad
D_{r,s}(y_0)\coloneqq   D_r(y_0)\times \torus^d_s\,.
\eeqano 

\item[\tiny $\bullet$] If  $\uno_d\coloneqq \diag(1)$ is the unit $(d\times d)$ matrix, we denote the standard symplectic matrix by 
$$\mathbb{J}\coloneqq  \begin{pmatrix}0 & -\uno_d\\
\uno_d & 0\end{pmatrix}\,.
$$ 

\item[\tiny $\bullet$]
For $y_0\in\real^d$, $\mathcal{A}_{r,s}(y_0)$ denotes the Banach space of real--analytic functions with bounded holomorphic extensions to $ D_{r,s}(y_0)$, with norm 
$$
\|\cdot\|_{r,s,y_0}\coloneqq \dst\sup_{ D_{r,s}(y_0)}|\cdot|\;.
$$
We also denote:
$$
\|\cdot\|_{r,y_0}\coloneqq \dst\sup_{ D_{r}(y_0)}|\cdot|\;,\quad \quad\|\cdot\|_{s}\coloneqq \dst\sup_{ \torus^d_s}|\cdot|\;.
$$
\item[\tiny $\bullet$] We equip 
$\cn\times\cn$   with the canonical symplectic form 
$$\varpi\coloneqq dy\wedge dx=dy_1\wedge dx_1+\cdots+dy_d\wedge dx_d\ ,
$$
and denote by $\phi_H^t$ the associated Hamiltonian flow governed by the Hamiltonian $H(y,x)$, $y,x\in\complex^d$, i.e., $z(t)\coloneqq \phi_H^t(y,x)$ is the solution of the Cauchy problem $\dot z= \mathbb{J} \nabla H(z)$, $z(0)=(y,x)$. 
\item[\tiny $\bullet$]
Given a linear operator $\mathcal{L}$ from the normed space $(V_1,\|\cdot\|_1)$ into the normed space  $(V_2,\|\cdot\|_2)$, its ``operator--norm'' is given by
\[\|\mathcal{L}\|\coloneqq \sup_{x\in V_1\setminus\{0\}}\frac{\|\mathcal{L}x\|_2}{\|x\|_1},\quad \mbox{so that}\quad \|\mathcal{L}x\|_2\le \|\mathcal{L}\|\, \|x\|_1\quad \mbox{for any}\quad x\in V_1.\]
\item[\tiny $\bullet$]
Given $\o\in \rn$, the directional derivative of a $C^1$ function $f$ with respect to $\o$ is given by
\[D_\o f\coloneqq \o\cdot  f_x=\dst\sum_{j=1}^d \o_j \dst f_{{x}_j}\,.\]
\item[\tiny $\bullet$]
If  $f$ is a (smooth or analytic)  function on $\torus^d$, its Fourier expansion is given by    \[f=\dst\sum_{k\in \zn}f_k \ex^{ik\cdot x}\,,
\qquad f_k\coloneqq\dst\frac{1}{(2\pi)^d}\dst\int_{\tn}f(x) \ex^{-ik\cdot x}\, d x\,,\]
(where, as usual, $\ex\coloneqq \exp(1)$ denotes the Neper number and $i$ the 
imaginary unit). We also set:
\[\average{f}\coloneqq f_0=\dst\frac{1}{(2\pi)^d}\dst\int_{\tn}f(x)\, d x\,,\qquad  ({\bf p}_N f)(x)\coloneqq\dst\sum_{|k|_1\leq N}f_k \ex^{ik\cdot x},\, N>0\,.
\]
${\bf p}_N$ being the Fourier projection onto the Fourier modes with $|k|_1\le N$; notice that $\langle\cdot\rangle= {\bf p}_0(\cdot)$.

\end{itemize}

\nl
We are ready to formulate a quantitative version of Arnold's Theorem\footnote{To avoid to introduce too many symbols, we use capital straight style for positive constants ($\mathsf{P}, 
\mathsf{K}, \mathsf{T}, \mathsf{C},...$), while, usually, capital normal style is used for functions or matrices ($K,P,H,T,...$). }:

\Giu
{\bf Theorem A} {\sl
Let $d \ge 2$; $\t\ge d-1$; $\a,r,\e>0$; $0<s_*<s\le 1$; $y_0\in\rn$;  $K,P\in \mathcal{A}_{r,s}(y_0)$; $H:=K+\e P$. Assume that
\beq{ArnoldCondv2prim}
\left\{ \begin{array}{l}
\o\coloneqq \dpr_y K(y_0)\in \D^\t_\a\,, \\ \ \\ 
 \det (\dpr^2_y K(y_0))\not= 0\;.
\end{array}
\right.
\eeq
Define:
$$
T\coloneqq \dpr^2_y K(y_0)^{-1},\ \mathsf{P}\coloneqq \|P\|_{r,s,y_0},\ \mathsf{K}\coloneqq \|\dpr^2_y K\|_{r,y_0},\  \mathsf{T}\coloneqq \|T\|\,,\ \torsion\coloneqq \mathsf{T}\mathsf{K}\,,
$$
and denote by $\m$ the rescaled smallness parameter:
\beq{epsilon}
\boxed{\m\coloneqq \mathsf{K}\mathsf{P} \frac{\vae\, }{\a^2}}\,.
\eeq
There exist constants $1<\mathsf{C}<\mathsf{C}_*$ depending only on $d$ and $\tau$, such that, if $a\coloneqq 6\t+3d+8$ and
\beq{smcondwhL}
\a\le \frac{r}{\mathsf{T}}\qquad\qquad \mbox{and}\qquad\qquad\m \le  \m_*\coloneqq\frac{(s-s_*)^{a}}{\mathsf{C}_*\;\torsion^4}\;,
\eeq
then, there exists a real--analytic  embedding 
$$\phi_*\colon x\in \torus^d_{s_*}\mapsto \phi_*(x)\coloneqq \phi_{\rm e}(y_0,x)+\big(v_*(x), u_*(x)\big) \in D_{r,s}(y_0)\,,
$$
where $\phi_{\rm e}$ is the trivial embedding 
\[ 
\phi_{\rm e}\colon x\in \tn \to (y_0,x),
\]
such that the  $d$--torus
\beq{KronTorArnv2}
\mathcal{T}_{\o,\vae}\coloneqq \phi_*\left(\tn\right)
\eeq
is a 
Lagrangian torus   satisfying
\beqno
\phi^t_{H}\circ \phi_*(x)=\phi_*(x+\o t) \,,\quad \forall \ x\in \torus^d_{s_*}\ ,\  \forall \ t\in\real\ .
\eeqno
Furthermore, 

\beq{est}
\max\Big\{ \| u_*\|_{s_*}\,, \frac{1}{2 \ex}\, \|\partial_x u_*\|_{s_*}\,,\, \frac{\mathsf{K}}{\a}\, \|v_*\|_{s_*}
\Big\} \le \frac{\mathsf{C}\  \torsion^3}{(s-s_*)^a} \ \m \le \frac{1}{4\ex}\,.
\eeq
}

\Giu
{\bf Remarks and addenda}\nobreak
\begin{itemize}
\item[(i)] $\torsion$ is a measure of the local ``torsion'' and is a number  greater than or equal to  one: 
\beq{Eta0Ge1}
\torsion\coloneqq \mathsf{T}\mathsf{K}\ge \mathsf{T}\|K_{yy}(y_0)\|\ge \|T\|\|K_{yy}(y_0)\|=\|T\|\|T^{-1}\|\ge 1\ .
\eeq

\item[(ii)] Notice that  the estimate on $v_*$ in \equ{est} implies that the maximal action oscillation of the torus $\mathcal{T}_{\o,\vae}$  is bounded by a constant times 
$\a \m$, which in view of \equ{epsilon}, is $\sim \e/\a$ as advertised in \equ{osc.int*}.

\item[(iii)] All numerical constants are explicitly ``computed'' during the proof. A complete list of them, including the definitions of $\mathsf{C}_*$ and $\mathsf{C}$,  is given in Appendix~\ref{AssumpExtArnolv2}.

\item[(iv)] The torus $\mathcal{T}_{\o,\vae}$ is {\sl Kolmogorov non--degenerate}. More precisely, $H$ can be put in Kolmogorov's normal form with non--degenerate quadratic part: there exists a 
symplectic transformation $\phi$ close to $\phi_{\rm e}$, for which
$$
H\circ \phi (y,x)= E+\o \cdot y + Q(y,x) \quad{\rm such\ that}\quad  \det \langle Q_{yy}(0,\cdot)\rangle\neq 0\,;
$$
for details, see Appendix~\ref{salamon}.

\item[(v)] The value of $\m_*$ in  \equ{smcondwhL} is not optimal. In Remark~\ref{rem:final} a better (still not optimal) value is given.

\item[(vi)] The dependence of the invariant torus $\mathcal{T}_{\o,\vae}$ on $\vae$ is analytic. More generally, if $H=H(y,x;z)$ is real--analytic also in $z\in V$, $V$ being some open set in $\complex^m$, and all the above norms are uniform in $z\in V$, then the invariant torus $\mathcal{T}_{\o,z}$ is real analytic in $V$. This is an obvious corollary of Weierstrass's theorem on uniform limits of holomorphic functions, in view 
of the uniformity of the limits in the proof.  

\end{itemize}

\section{Proof}
\subsection{Arnold's scheme: the basic step}\label{Astep}
The next Lemma describes Arnold's basic KAM step, on which Arnold's scheme is based. Its quantitative formulation involves a few constants, which are defined as follows:
\beqano%* {costants1}
\n &\coloneqq& \t+1\;,\quad
\mathsf{C}_0 \coloneqq 4\sqrt{2}\left(\frac{3}{2}\right)^{2\n+d}\dst\int_{\rn} \left( |y|_1^{\n}+|y|_1^{2\n}\right)\ex^{-|y|_1}dy\;,\\
\mathsf{C}_1 &\coloneqq& 2\left(\frac{3}{2}\right)^{\n+d}\dst\int_{\rn} |y|_1^{\n}\ex^{-|y|_1}dy\;,\\
\mathsf{C}_2 &\coloneqq& 2^{3d}d\;,\quad
\mathsf{C}_3 \coloneqq \left(d^2\mathsf{C}_1^2+6d\mathsf{C}_1 +\mathsf{C}_2\right)\sqrt{2}\;,\quad
\mathsf{C}_4 \coloneqq \max\left\{6d^2\mathsf{C}_0,\,\mathsf{C}_3\right\}\;.
\eeqano

\lem{lem:1bisv2}
Let\footnote{$K$ and $P$  stand, here,  for  generic real analytic Hamiltonians which, later on, will respectively play the roles of $K_j$ and $P_j$,  and $\mathsf{y},\,r$, the roles of $y_j,\,r_j$ in the iterative step.}
 $r>0,\,0<2\s<s\leq 1$, $\mathsf{y}\in\real^d$,  $K,P\in \mathcal{A}_{r,s}(\mathsf{y})$ and consider the Hamiltonian parametrised by 
%* $\vae\ge 0$
$\vae> 0$
$$
H(y,x;\vae)\coloneqq K(y)+\vae P(y,x)\,.
$$
Assume that 
$$
\det K_{yy}(\mathsf{y})\neq 0\ ,\qquad \quad  \o\coloneqq K_{y}(\mathsf{y})\in \D^\t_\a\ ,
$$
and let 
$ \mathsf{K}$, $ \mathsf{T}$ and $\mathsf{P}$ be positive numbers such that
\beq{RecHypArnv2}
\|K_{yy}\|_{r,\mathsf{y}}\le \mathsf{K}\;, \qquad\ \; \|T\|\le \mathsf{T}\;,\qquad
\|P\|_{r,s,\mathsf{y}}\le \mathsf{P} \;,
\eeq
where $T\coloneqq K_{yy}(\mathsf{y})^{-1}$.
\\
Now, let $\l, \check{r}, \bar r$ be positive number such that:
\beq{lamsup1}
\l\ge \log\Big(\s^{2\n+d}\frac{{\a}^2}{\vae\,{\mathsf{P}}\mathsf{K}}\Big)\;, \quad 
\check{r}\le \frac{5}{24d}\frac{r}{\mathsf{T} \mathsf{K}}\;,\quad 
\bar{r}\le
\dst\min\left\{\check{r}   \,,\,   \frac{\a}{2d\mathsf{K}\k^{\t+1}}  \right\}\ ,
\eeq
where
$$
\k\coloneqq \frac{4\l}{\s}\ .
$$
Finally, define 
$$
\mathsf{L}\coloneqq \mathsf{P}\dst\max\Big\{\frac{40d\mathsf{T}^2\mathsf{K}  }{r^2}\s^{-(\n+d)}\,,\,\frac{\mathsf{C}_4}{\sqrt{2}} \max\big\{1,\frac{\a}{r\mathsf{K}}\big\}\frac{ \mathsf{K}}{\a^2}\s^{-2(\n+d)}\Big\}\ ,\quad
\bar{s}\coloneqq s-\frac{2}{3}\s,  \quad s'\coloneqq s-\s \ .
$$
Then, if
\beq{cond1Bisv2}
{\vae\, }{\mathsf{L}}\le \frac{\sigma}{3}\;,
\eeq
there exist $\mathsf{y}'\in\rn$ and a symplectic change of coordinates
\beq{phiokBis0v2}
\phi'=\id + \vae \tilde \phi\colon D_{\bar{r}/2,s'}(\mathsf{y}')\to D_{2r/3, \bar{s}}(\mathsf{y}),
\eeq
%* generated by a function $g\in \mathcal{A}_{\bar r,\bar s}(\mathsf{y})$ 
such that
\beq{HPhiH'}
\left\{
\begin{aligned}
& H\circ \phi'\eqqcolon H'\eqqcolon K'+\vae^2 P'\ ,\\
& \dpr_{y'} K'(\mathsf{y}')=\o,\quad \det \dpr^2_{y'} K'(\mathsf{y}')\neq 0\,,
\end{aligned}
\right.
\eeq
where 
$$K'\coloneqq K+\vae\wt K\coloneqq K+\vae \average{P(y',\cdot)}\;.
$$
Moreover, letting 
$$\left(\dpr^2_{y'} K'(\mathsf{y}')\right)^{-1}\eqqcolon T+\vae\;\wt T\;,
$$
the following estimates hold:
\beqa{convEstv2}
&&\|\dpr_{y'}^2\wt K\|_{r/2,\mathsf{y}}\le 
\mathsf{K}\mathsf{L}\,,\quad |\mathsf{y}'-\mathsf{y}| 
\le \frac{8\vae\, \mathsf{T} \mathsf{P}}{r}\,,\quad
\|\wt T\|\le \mathsf{T}\mathsf{L}
\,, \\
&&  \max\{\|\dpr_x\pi_2\tilde \phi\|_{s'},\|\mathsf{W}\,\tilde \phi\|_{\bar{r}/2,s',\mathsf{y}'}\}\le d^{-2}\s^{d-1}{\mathsf{L}}\,,\quad
\|P'\|_{\bar{r}/2, s',\mathsf{y}'}\le  \mathsf{L}\mathsf{P}\;,
\eeqa
where 
$$
\mathsf{W}\coloneqq \begin{pmatrix}
\max\{\frac{\mathsf{K}}{{\a}}\;,\frac{1}r\}\;\uno_d & 0\\ \ \\
0			& \uno_d 
\end{pmatrix}.
$$
\elem
Observe that 
$$\s^{-2(\n+d)}\vae\;\mathsf{P}\mathsf{K}/\a^2\le (\sqrt{2}/\mathsf{C}_4)\;\vae\;\mathsf{L}\,,
$$
so that \equ{cond1Bisv2} implies 
\beqno
\frac{\vae\,{\mathsf{P}}\mathsf{K}}{{\a}^2}<\frac{\s^{2\n+d}}{\ex}\;,%* warning: this assumption has to be checked
\eeqno
which, in particular, implies that $\l>1$ and $\k>4$.

%*
%For later use, observe that 
%$$\mathsf{L}> \s^{-d+1}\ovl{\mathsf{L}}> \ovl{\mathsf{L}}$$ 
%and, also, that:  
%$$\frac{40d\mathsf{T}^2\mathsf{K} }{r^2}\s^{-(\n+d)}>\frac{4\mathsf{T} }{r^2}\ge\frac{4}{\mathsf{K} r^2}\;.
%$$
%

\proof

\nl
{\bf Step 1: Construction of  Arnold's transformation }  

\nl
We seek  a near--identity symplectic transformation 
\[\phi'\colon D_{r_1,s_1}(\mathsf{y}')\to D_{r,s}(\mathsf{y}),\]
with $D_{r_1,s_1}(\mathsf{y}')\subset D_{r,s}(\mathsf{y})$,   generated by a generating function\footnote{Following the classical approach of Arnold, we use generating functions to construct symplectic transformations. Of course one could also use the equivalent method of time--one Hamiltonian flows (or Lie series).}
of the form $y'\cdot x+\vae g(y',x)$, so that
\beq{ArnTraKamv2}
\phi'\colon \left\{\begin{aligned}
y  &=y'+\vae g_x(y',x)\\
x' &=x+\vae g_{y'}(y',x)\, ,
\end{aligned}
\right.
\eeq
such that
\beq{ArnH1v2}
\left\{
\begin{aligned}
& H'\coloneqq H\circ \phi'=K'+\vae^2 P'\ ,\\
& \dpr_{y'} K'(\mathsf{y}')=\o,\quad \det \dpr^2_{y'} K'(\mathsf{y}')\neq 0\,.
\end{aligned}
\right.
\eeq
By Taylor's formula, we get\footnote{Recall (\S 2) that $\average{\cdot}$ stands for the average over $\tn$ and that ${\bf p}_N$ is the Fourier projection onto modes with $|k|_1\le N$.}
\beq{Arneq11v2}
\begin{aligned}
H(y'+\vae g_x(y',x),x)=&K(y')+\vae \wt K(y') +\vae \left[K'(y')\cdot g_x +{\bf p}_{\k} P(y',\cdot)-\wt K(y') \right]+\\
						&+\vae^2 \left( P^\ppu+P^\ppd+ P^\ppt\right)(y',x) \\
			= & K'(y')+\vae \left[K'(y')\cdot g_x +{\bf p}_{\k} P(y',\cdot)-\wt K(y') \right]+ \vae^2 P_+(y',x),
\end{aligned}
\eeq
with $\k>0$, %* warning: we don't need natural numbers
 which will be chosen large enough so that $P^\ppt=O(\vae)$ and 
\beq{ArnDefPsv2}
\left\{
\begin{aligned}
P_+\ &\coloneqq P^\ppu+P^\ppd+ P^\ppt\\
P^\ppu &\coloneqq \su{\vae^2}\left[K(y'+\vae g_x)-K(y')-\vae K_y(y')\cdot g_x \right]=\dst\int^1_0(1-t)K_{yy}(\vae t g_x)\cdot g_x\cdot g_x dt\\
P^\ppd &\coloneqq \su\vae \left[P(y'+\vae g_x,x)-P(y',x)\right]=\dst\int_0^1P_y(y'+\vae t g_x,x)\cdot g_x dt\\
P^\ppt &\coloneqq \su\vae \left[ P(y',x)-{\bf p}_{\k} P(y',\cdot)\right]=\su\vae \dst\sum_{|n|_1>\k} P_n(y')\ex^{in\cdot x}\; .
\end{aligned}
\right.
\eeq
By the non--degeneracy condition $\det K_{yy}(\mathsf{y})\neq 0$, for $\vae$ small enough (to be made precised below), $\det\dpr_{y'}^2 K'(\mathsf{y})\neq0$
 and, therefore, by the standard Inverse Function Theorem (see, e.g., Lemma~\ref{IFTLem}), there exists a unique $\mathsf{y}'\in D_r(\mathsf{y})$ such that the second part of \eqref{ArnH1v2} holds. 
In view of \eqref{Arneq11v2}, in order to get the first part of \eqref{ArnH1v2}, we need to find $g$ such that  $K_y(y')\cdot g_x +{\bf p}_{\k} P(y',\cdot)-\wt K(y')$ vanishes; such a $g$ is indeed given by
 \beq{HomEqArnv2}
 g\coloneqq \dst\sum_{0<|n|_1\leq \k} \frac{-P_n(y')}{iK_y(y')\cdot n}\ex^{in\cdot x},
 \eeq
provided that 
\beq{CondHomEqArnv2}
K_y(y')\cdot n\neq 0, \quad \forall\; 0<|n|_1\leq \k,\quad \forall\; y'\in D_{r_1}(\mathsf{y}')\quad  \left(\subset D_{r}(\mathsf{y})\right).
\eeq
But, in fact, since $K_y(\mathsf{y})$ is rationally independent, then, given any $\k>0$, % as above
 there exists $\bar{r}\leq r$ such that
\beq{CondHomEqArnBisv2}
K_y(y')\cdot n\neq0,\quad \forall\; 0<|n|_1\leq \k, \quad\forall\; y'\in D_{\bar{r}}(\mathsf{y}).
\eeq
The last step is to invert the function $x\mapsto x+\vae g_{y'}(y',x)$ in order to define $P'$. By the Inverse Function Theorem, for $\vae$ small enough, the map $x\mapsto x+\vae g_{y'}(y',x)$ admits a real--analytic inverse of the form
\beq{InvComp2Fiv2}
\f_\vae(y',x')= x'+\vae \wt \f_\vae(y',x'),
\eeq
so that the Arnold's symplectic transformation is given by
\beq{ArnTrans0v2}
\phi'\colon (y',x')\mapsto \left\{
\begin{aligned}
y &= y'+\vae g_x(y',\f_\vae(y',x'))\\
x &= \f_\vae(y',x')= x'+\vae \wt \f_\vae(y',x') .
\end{aligned}
\right.
\eeq
Hence, \eqref{ArnH1v2} holds with
\beq{DefP1Arv2}
P'(y',x')\coloneqq P_+(y', \f_\vae(y',x')).
\eeq
{\bf Step 2: Quantitative estimates}\\
First of all, notice that from the definitions of $\bar r$ and $\check r$ it follows that
\beq{rrbarAsv2}
\bar{r}\le \check r\le \frac{5r}{24d}<\frac{r}{2}\;.
\eeq
\noi
We begin by extending the ``Diophantine condition w.r.t. $K_y$'' uniformly to $D_{\bar{r}}(\mathsf{y})$ up to the order $\k$. Indeed, by the Mean Value Inequality and $K_y(\mathsf{y})=\o\in\D^\t_\a$,  we get, for any $0<|n|_1\leq \k$ and any $y'\in D_{\bar{r}}(\mathsf{y})$,
\begin{align}
|K_y(y')\cdot n|&=|\o\cdot n +(K_y(y')-K_y(\mathsf{y}))\cdot n|\geq |\o\cdot n|\left(1-d\frac{\|K_{yy}\|_{\bar{r},\mathsf{y}}}{|\o\cdot n|}|n|_1\bar{r}\right) \nonumber\\
         &\geq \frac{\a}{|n|_1^\t}\left(1-\frac{d\mathsf{K}}{\a }|n|_1^{\t+1}\bar{r} \right)\geq \frac{\a}{|n|_1^\t}\left(1-\frac{d\mathsf{K}}{\a }\k^{\t+1}\bar{r} \right)\ge \frac{\a}{2|n|_1^\t},\label{ArnExtDiopCondv2}
\end{align}
so that, by Fourier estimates (Lemma~\ref{Cau}--{\bf (ii)}), we have
\begin{align*}
\|g_x\|_{\bar{r},\bar{s},\mathsf{y}} &\overset{def}{=}\dst\sup_{D_{\bar{r},\bar{s}}(\mathsf{y})}\left|\dst\sum_{0<|n|_1\leq \k}\frac{nP_n(y')}{K_y(y')\cdot n}\ex^{in\cdot x} \right|\leq \dst\sum_{0<|n|_1\leq \k}\frac{\|P_n\|_{\bar{r},\bar{s}, \mathsf{y}}}{|K_y(y')\cdot n|}|n|_1\ex^{\left(s-\frac{2}{3}\s\right)|n|_1}\\
   &\leq \dst\sum_{0<|n|_1\leq \k} \mathsf{P}\ex^{-s|n|_1}\frac{2|n|_1^{\n}}{\a}\ex^{\left(s-\frac{2}{3}\s\right)|n|_1}\leq \frac{2\mathsf{P}}{\a}\dst\sum_{n\in\zn} |n|_1^{\n}\ex^{-\frac{2}{3}\s|n|_1}\\
   &\leq \frac{2\mathsf{P}}{\a}\dst\int_{\rn} |y|_1^{\n}\ex^{-\frac{2}{3}\s|y|_1}dy
 = \left(\frac{3}{2\s}\right)^{\n+d}\frac{2\mathsf{P}}{\a}\dst\int_{\rn} |y|_1^{\n}\ex^{-|y|_1}dy
   = \mathsf{C}_1 \frac{\mathsf{P}}{\a} \s^{-(\n+d)}\,,
\end{align*}
\begin{align*}
\|\dpr_{y'}g\|_{\bar{r},\bar{s},\mathsf{y}} &\overset{def}{=}\dst\sup_{D_{\bar{r},\bar{s}}(\mathsf{y})}\left|\dst\sum_{0<|n|_1\leq \k}\left(\frac{ \dpr_yP_n(y')}{K_y(y')\cdot n}-P_n(y')\frac{ K_{yy}(y')n}{(K_y(y')\cdot n)^2}\right)\ex^{in\cdot x} \right|\\
   &\leq \dst\sum_{0<|n|_1\leq \k}\dst\sup_{D_{\bar{r}}(\mathsf{y})}\left(\frac{\|(P_y)_n\|_{\bar{r},s, \mathsf{y}}}{|K_y(y')\cdot n|}+\|P_n\|_{r,s, \mathsf{y}}\frac{\|K_{yy}\|_{r,\mathsf{y}}|n|_1}{|K_y(y')\cdot n|^2}\right)\ex^{\left(s-\frac{2}{3}\s\right)|n|_1}\\
   &\stackrel{\equ{RecHypArnv2}+\equ{ArnExtDiopCondv2}}{\le} \dst\sum_{0<|n|_1\leq \k}\left( \frac{\mathsf{P}}{r-\bar{r}}\ex^{-s|n|_1}\frac{2|n|_1^{\t}}{\a}+\mathsf{P}\ex^{-s|n|_1}\mathsf{K}|n|_1\left(\frac{2|n|_1^{\t}}{\a}\right)^2\right)\ex^{\left(s-\frac{2}{3}\s\right)|n|_1}\\
      &\leby{rrbarAsv2} \frac{4\mathsf{P}}{\a^2 r}\dst\sum_{0<|n|_1\leq \k}\left( |n|_1^{\t}\a +r\mathsf{K}|n|_1^{2\t+1}\right)\ex^{-\frac{2}{3}\s|n|_1}\\
   &\le \max\left\{\a,r\mathsf{K}\right\}\frac{4\mathsf{P}}{\a^2 r}\dst\sum_{0<|n|_1\leq \k}\left( |n|_1^{\t}+|n|_1^{2\t+1}\right)\ex^{-\frac{2}{3}\s|n|_1}\\
   &\leq \max\left\{1,\frac{\a}{r\mathsf{K}}\right\}\frac{4\mathsf{P} \mathsf{K}}{\a^2 }\dst\int_{\rn} \left( |y|_1^{\t}+|y|_1^{2\t+1}\right)\ex^{-\frac{2}{3}\s|y|_1}dy \\
   &= \left(\frac{3}{2\s}\right)^{2\t+d+1}\max\left\{1,\frac{\a}{r\mathsf{K}}\right\}\frac{4\mathsf{P} \mathsf{K}}{\a^2 }\dst\int_{\rn} \left( |y|_1^{\t}+|y|_1^{2\t+1}\right)\ex^{-|y|_1}dy\\
   &\le \frac{\mathsf{C}_0}{\sqrt{2}} \max\left\{1,\frac{\a}{r\mathsf{K}}\right\}\frac{\mathsf{P} \mathsf{K}}{\a^2 }\s^{-(2\t+d+1)}
   < \ovl{\mathsf{L}} \;,
\end{align*}
where
$$
\ovl{\mathsf{L}}\coloneqq 6\ \frac{\mathsf{C}_0}{\sqrt{2}} \max\left\{1,\frac{\a}{r\mathsf{K}}\right\}\frac{\mathsf{P} \mathsf{K}}{\a^2}\s^{-(2\n+d+1)}\ .
$$
Analogously,
\begin{align*}
\|\dpr^2_{y'x}g\|_{\bar{r},\bar{s},\mathsf{y}} 
   \le \frac{\mathsf{C}_0}{\sqrt{2}} \max\left\{1,\frac{\a}{r\mathsf{K}}\right\}\frac{\mathsf{P} \mathsf{K}}{\a^2 }\s^{-(2\n+d)}
   \le\ovl{\mathsf{L}}\;,
\end{align*}
and, by Cauchy's estimate (Lemma~\ref{Cau}--{\bf (i)}) we get
\beq{gy'xx}
\|\dpr^3_{y'xx}g\|_{\bar{r},s'',\mathsf{y}} 
   \le \frac{6\mathsf{C}_0}{\sqrt{2}} \max\left\{1,\frac{\a}{r\mathsf{K}}\right\}\frac{\mathsf{P} \mathsf{K}}{\a^2 }\s^{-(2\n+d+1)}
   =\ovl{\mathsf{L}}\;,
\eeq
where
$$
s''\coloneqq s-\frac{5}{6}\s
\quad {\rm and}\quad
\|\dpr^3_{y'xx}g\|_{\bar{r},s'',\mathsf{y}}\coloneqq \sup_{D_{\bar{r},s''}(\mathsf{y})}\max\{|\dpr^3_{y'_ix_jx_k}g|\ :\ i,j,k=1,\cdots,d\}\;.
$$ 
Also,
\[\|\wt K_y\|_{r/2,\mathsf{y}}=\| \langle P_y\rangle \|_{r/2,\mathsf{y}}\leq \|P_y\|_{r/2,\bar{s}, \mathsf{y}}\leq  \frac{\mathsf{P}}{r-\frac{r}{2}}\leq \frac{2\mathsf{P}}{r} \;,\]
\[\|\dpr_{y'}^2\wt K\|_{r/2,\mathsf{y}}=\| \langle P_{yy}\rangle\|_{r/2,\mathsf{y}}\leq \|P_{yy}\|_{r/2,\bar{s}, \mathsf{y}}\leq  \frac{\mathsf{P}}{(r-\frac{r}{2})^2}\leq \frac{4\mathsf{P}}{r^2}\le \mathsf{K}\mathsf{L} 
\;.\]
Next, we prove the existence and uniqueness of $\mathsf{y}'$ in \eqref{ArnH1v2}. 
Let $U_{\vae}\coloneqq \{\eta\in\complex: |\eta|< 2\vae\, \}$ and consider the map:
\begin{align*}
F\colon D_{\check{r}}(\mathsf{y})\times U_{\vae} &\longrightarrow \qquad \cn\\
		(y,\eta)\quad &\longmapsto K_y(y)+\eta \wt K_{y'}(y)-K_y(\mathsf{y})\;.
\end{align*} 
Then
\begin{itemize}
\item $F(\mathsf{y},0)=0,\quad F_y(\mathsf{y},0)^{-1}=K_{yy}(\mathsf{y})^{-1}=T$.
\item For any $(y,\eta)\in D_{\check{r}}(\mathsf{y})\times U_{\vae}$,
\begin{align*}
\|\uno_d-TF_y(y,\eta)\|&\leq \|\uno_d-TK_{yy}\|+|\eta|\;\|T\|\;\|\dpr_{y'}^2\wt K\|_{r/2,\mathsf{y}}\\
	  &\leq d\|T\|\|K_{yyy}\|_{\check{r},\mathsf{y}}\check{r}+ 2\vae\, \mathsf{T}\frac{4\mathsf{P}}{r^2}\\
      &\leq d\mathsf{T} \mathsf{K}\frac{\check{r}}{r-\check{r}}+8\mathsf{T}\frac{\vae\, \mathsf{P}}{ r^2}
      \leby{rrbarAsv2}d\mathsf{T} \mathsf{K} \frac{2\check{r}}{ r}+\vae\, \frac{8\mathsf{T} \mathsf{P}}{r^2}\\
      &\le 2d\mathsf{T} \mathsf{K}\frac{\bar{r}}{ r}+\su2{\vae\, }\mathsf{L}\\
      &\overset{\equ{rrbarAsv2}+\equ{cond1Bisv2}}{\leq}\frac{5}{12}+\frac{\s}{6}
      \le \frac{5}{12}+\su{12}=\su2\;.
\end{align*}
\item Recalling $\s\le\su2$, we have
\begin{align}
2\|T\|\|F(\mathsf{y},\cdot)\|_{2\vae\, ,0}&=2\|T\|\dst\sup_{U_{\vae}}|\eta \wt K_{y'}(\mathsf{y})|
		\leq 2\mathsf{T} \frac{4\vae\, \mathsf{P}}{r}
		\le \frac{5\cdot 2^{\n+d}}{8d}\frac{r}{\mathsf{T}\mathsf{K}}\s^{\n+d}{\vae\, }\mathsf{L}\nonumber\\
		&= 3\cdot 2^d\;(2\s)^{\n}\;\check r\;\s^{d}{\vae\, }\mathsf{L} 
		\le 3\cdot 2^d\;\check r\;\s^{d}{\vae\, }\mathsf{L}\label{distyy1I}\\
		&\leby{cond1Bisv2} 3\;\check r\;(2\s)^{d}\;\frac{\s}{3} \le \frac{\check{r}}{2}\;.\nonumber
\end{align}
\end{itemize}
Therefore, we can apply the Inverse Function Theorem (Lemma~\ref{IFTLem}). Hence, there exists a function $g\colon U_{\vae}\to D_{\check{r}}(\mathsf{y})$ such that its graph coincides with $F^{-1}(\{0\})$. In particular, $\mathsf{y}'\coloneqq g(\vae)$ is the unique $y\in D_{\check{r}}(\mathsf{y})$ satisfying $0=F(y,\vae)=\dpr_y K'(y)-\o$, i.e.,  the second part of \eqref{ArnH1v2}. Moreover, 
\beq{EcarY1Y0v2}
|\mathsf{y}'-\mathsf{y}|\leq 2\|T\|\|F(\mathsf{y},\cdot)\|_{2\vae\, ,0}\leq \frac{8\vae\, \mathsf{T} \mathsf{P}}{r}\leby{distyy1I} 3\cdot 2^d\;\check r\;\s^{d}{\vae\, }\mathsf{L}\leq \frac{\check{r}}{2}\;,
\eeq
so that
\beq{NextSetArnv2}
D_{\frac{\check{r}}{2}}(\mathsf{y}')\subset D_{\check{r}}(\mathsf{y}).
\eeq
Next, we prove that $\dpr^2_y K'(\mathsf{y}')$ is invertible. Indeed, by Taylor' formula, we have
\begin{align*}
\dpr^2_y K'(\mathsf{y}')&= K_{yy}(\mathsf{y})+ \dst\int_0^1 K_{yyy}(\mathsf{y}+t\vae \wt y)\cdot\vae\wt y dt+\vae \wt K_{yy}(\mathsf{y}')\\
           &= T^{-1}\left(\uno_d+\vae T\left(\dst\int_0^1 K_{yyy}(\mathsf{y}+t\vae \wt y)\cdot\wt y dt+ \wt K_{yy}(\mathsf{y}')\right)\right)\\
           &\eqqcolon T^{-1}(\uno_d+\vae A),
\end{align*}
and, by Cauchy's estimate, 
\begin{align*}
\vae\, \|A\|&\leq \|T\|\left(d\|K_{yyy}\|_{r/2,\mathsf{y}}\vae\, |\mathsf{y}'-\mathsf{y}|+ \vae\, \|\dpr_{y'}^2\wt K\|_{r/2,\mathsf{y}}\right)\\
     &\leq \|T\|\left(\frac{d\|K_{yy}\|_{r,\mathsf{y}}}{r-\frac{r}{2}}\vae\, |\mathsf{y}'-\mathsf{y}|+\vae\, \|\wt K_{yy}\|_{r/2,\mathsf{y}}\right)\\
	 &\leby{EcarY1Y0v2} \mathsf{T}\left(\frac{2d\mathsf{K}}{r}\frac{8\vae\, \mathsf{T} \mathsf{P}}{r}+\frac{4\vae\, \mathsf{P}}{r^2} \right)
	 \leq \frac{4\vae\, \mathsf{T}\mathsf{P}}{r^2}(4d\mathsf{T}\mathsf{K}+1)\\
	 &\leq\frac{20d\vae\, \mathsf{T}^2\mathsf{K} \mathsf{P}}{r^2}
	 \le \su 2\vae\, \mathsf{L}
	 \leby{cond1Bisv2}\frac{\s}{6}
	 \le\su2.
\end{align*}
Hence $\dpr_{y'}^2 K'(\mathsf{y}')$ is invertible with
\[\dpr_{y'}^2 K'(\mathsf{y}')^{-1}=(\uno_d+\vae A)^{-1}T=T+\dst\sum_{k\geq 1}(-\vae)^k A^k T\eqqcolon T+\vae \wt T,\]
and
\[\vae\, \|\wt T\|\leq \vae\, \frac{\|A\|}{1-\vae\, \|A\|}\|T\|\leq 2\vae\, \|A\| \|T\|
\le \vae\, \mathsf{L}\mathsf{T}
\le 2\frac{\s}{6}\mathsf{T}
= \mathsf{T}\frac{\s}{3}\,.\]
Next, we prove estimate on $P_+$. We have,
\[\vae\, \|g_x\|_{\bar{r},\bar{s},\mathsf{y}}\leq \vae\, \mathsf{C}_1 \frac{\mathsf{P}}{\a} \s^{-(\t+d+1)} \le \vae\, \frac{r}{3}\mathsf{L}\leby{cond1Bisv2}\frac{r}{3}\frac{\s}{3}\le \frac{r}{3}\]
so that, for any $(y',x)\in D_{\bar{r},\bar{s}}(\mathsf{y})$,
\[ |y'+\vae g_x(y',x)-\mathsf{y}|\leq \bar{r}+\frac{r}{3}< \frac{r}{8d}+\frac{r}{3}<\frac{2r}{3}<r\,,\]
and thus
\begin{align*}
\|P^\ppu\|_{\bar{r},\bar{s},\mathsf{y}}&\leq d^2 \|K_{yy}\|_{r,\mathsf{y}}\|g_x\|_{\bar{r},\bar{s},\mathsf{y}}^2\leq d^2 \mathsf{K}\left( \mathsf{C}_1 \frac{\mathsf{P}}{\a} \s^{-(\n+d)}\right)^2
   =d^2\mathsf{C}_1^2 \frac{\mathsf{K}\mathsf{P}^2}{\a^2} \s^{-2(\n+d)}, 
 \end{align*}
\begin{align*}
\|P^\ppd\|_{\bar{r},\bar{s},\mathsf{y}}&\leq d\|P_y\|_{\frac{5r}{6},\bar{s},\mathsf{y}}\|g_x\|_{\bar{r},\bar{s},\mathsf{y}}\leq d\frac{6\mathsf{P}}{r}\mathsf{C}_1 \frac{\mathsf{P}}{\a} \s^{-(\n+d)}
     = 6d\mathsf{C}_1 \frac{\mathsf{P}^2}{\a r}\s^{-(\n+d)}
\end{align*}
and by Fourier estimates (Lemma~\ref{Cau}--{\bf (ii)}), we have,
\begin{align*}
\vae\, \|P^\ppt\|_{\bar{r},s-\frac{\s}{2},\mathsf{y}}&\leq \dst\sum_{|n|_1>\k}\|P_n\|_{\bar{r},\mathsf{y}}\ex^{(s-\frac{\s}{2})|n|_1}\leq \mathsf{P}\dst\sum_{|n|_1>\k}\ex^{-\frac{\s |n|_1}{2}}\\
  &\leq \mathsf{P}\ex^{-\frac{ \k\s}{4}}\dst\sum_{|n|_1>\k}\ex^{-\frac{\s |n|_1}{4}}\leq \mathsf{P}\ex^{-\frac{ \k\s}{4}}\dst\sum_{|n|_1>0}\ex^{-\frac{\s |n|_1}{4}}\\
  &= \mathsf{P}\ex^{-\frac{ \k\s}{4}} \left(\left(\dst\sum_{k\in \integer}\ex^{-\frac{\s |k|}{4}}\right)^d-1\right)=\mathsf{P}\ex^{-\frac{ \k\s}{4}}\left(\left(1+\frac{2\ex^{-\frac{\s }{4}}}{1-\ex^{-\frac{\s }{4}}} \right)^d-1\right)\\
  &= \mathsf{P}\ex^{-\frac{ \k\s}{4}}\left(\left(1+\frac{2}{\ex^{\frac{\s }{4}}-1} \right)^d-1\right)\leq \mathsf{P}\ex^{-\frac{ \k\s}{4}}\left(\left(1+\frac{2}{\frac{\s }{4}} \right)^d-1\right)\\
  &\leq \s^{-d} \mathsf{P}\ex^{-\frac{ \k\s}{4}}\left(\left(\s +8 \right)^d-\s^d\right)\leq d 8^{d}\s^{-d} \mathsf{P}\ex^{-\frac{ \k\s}{4}}\\
  &= \mathsf{C}_2\s^{-d} \mathsf{P}\ex^{-\l}
  \leby{lamsup1} \mathsf{C}_2\s^{-d} \mathsf{P}\s^{-(2\n+d)}\frac{\vae\, {\mathsf{P}}\mathsf{K}}{{\a}^2}
  = \mathsf{C}_2 \mathsf{P}\frac{\vae\, {\mathsf{P}}\mathsf{K}}{{\a}^2}\s^{-2(\n+d)}\,.
\end{align*}
Hence, 
\begin{align*}
\|P_+\|_{\bar{r},\bar{s},\mathsf{y}}&\leq \|P^\ppu\|_{\bar{r},\bar{s},\mathsf{y}}+\|P^\ppd\|_{\bar{r},\bar{s},\mathsf{y}}+\|P^\ppt\|_{\bar{r},\bar{s},\mathsf{y}}\\
  &\leq d^2\mathsf{C}_1^2 \frac{\mathsf{K}\mathsf{P}^2}{\a^2} \s^{-2(\n+d)}+6d\mathsf{C}_1 \frac{\mathsf{P}^2}{\a r}\s^{-(\n+d)}+\mathsf{C}_2 \mathsf{P}\frac{\vae\, {\mathsf{P}}\mathsf{K}}{{\a}^2}\s^{-2(\n+d)}\\
  &= \left(d^2\mathsf{C}_1^2 r\mathsf{K}+6d\mathsf{C}_1 \a \s^{\n+d}+\mathsf{C}_2 r\mathsf{K}\right)\frac{\mathsf{P}^2}{\a^2 r}\s^{-2(\t+d+1)}\\
  &\le \left(d^2\mathsf{C}_1^2+6d\mathsf{C}_1 +\mathsf{C}_2\right)\max\left\{\a,r\mathsf{K}\right\}\frac{\mathsf{P}^2}{\a^2 r}\s^{-2(\t+d+1)}\\
  &\le 
  %* \leby{RecHypArnv2} 
  \frac{\mathsf{C}_3}{\sqrt{2}} \max\left\{1,\frac{\a}{r\mathsf{K}} \right\}\frac{\mathsf{P}^2\mathsf{K}}{\a^2 }\s^{-2(\n+d)}
  \le \mathsf{L}\mathsf{P}\;.
\end{align*}
Finally, we prove that, given $y'\in D_{\bar{r}}(\mathsf{y})$, the function $\psi_\vae(x)=x+\vae g_{y'}(y',x)$ has an analytic inverse\footnote{Observe that $\psi_\vae(id+\vae u)=id$ is equivalent to $u=-g_{y'}(y',id+\vae u)$, i.e., $u$ is a fixed--point of the map $u\mapsto-g_{y'}(y',id+\vae u)$.}.
Consider the Banach's space
$$
\mathcal{B}\coloneqq \bigg\{u\in C^1(\torus^d_{s'},\cn): \|u\|_{s',1}\coloneqq\max\big\{\|u\|_{s'}\;, \|\dpr_x u\|_{s'}\big\}\le \overline{\mathsf{L}} \bigg\}\;. 
$$
 For any $u\in\mathcal{B}$ and any $x'\in \torus^d_{s'}$, we have $\Im(x'+\vae u(x'))\le s'+\vae\, \|u\|_{s'}\le s'+\vae\, \overline{\mathsf{L}}\leby{cond1Bisv2} s'+\s/6=s''.$ Hence, the functional $f\colon \mathcal{B}\ni u\mapsto -g_{y'}(y',\id+\vae u)$ is well--defined and smooth. Moreover, for any $u\in \mathcal{B},$
$$
\|f(u)\|_{s'}\le \|g_{y'}\|_{W''}\le\overline{\mathsf{L}},\quad \|\dpr_x (f(u))\|_{s'}\le \|g_{y'x}\|_{W''}\cdot|\vae|\|\dpr_x u\|_{s'}\le\overline{\mathsf{L}}\cdot |\vae| \overline{\mathsf{L}}\leby{cond1Bisv2} \overline{\mathsf{L}}\cdot \frac{\s}{6}<\ovl{\mathsf{L}}. 
$$
Thus, $f\colon \mathcal{B}\to\mathcal{B}$. Furthermore, for any $u_1,u_2\in\mathcal{B}$, 
$$\|f(u_1)-f(u_2)\|_{s',1}\le (1+d^2\vae\, \overline{\mathsf{L}}) \vae\, \overline{\mathsf{L}}\cdot\|u_1-u_2\|_{s',1}\leby{cond1Bisv2}2\frac{\s}{3d^2}\cdot \|u_1-u_2\|_{s',1}<\su{2} \|u_1-u_2\|_{s',1},
$$
Hence, $f$ is a contraction. Therefore, by the Banach--Caccioppoli fixed--point Theorem, $f$ has a unique fixed--point $\wt{\f}_\vae\in \mathcal{B}$; $\wt{\f}_\vae$ is obtained as the uniform limit $\dst\lim_n f^n(0)$ (as $0\in \mathcal{B}$). Thus, as $f^0= f$ is real--analytic on $D_{\bar{r}}(\mathsf{y})\times\torus^d_{s'}$, by Weierstrass's Theorem on the uniform convergence of analytic functions, $\wt{\f}_\vae$ is real--analytic on $D_{\bar{r},s'}(\mathsf{y})$. The rest of the claims on $\phi'$ and $P'$ are then obvious.
\qed

\subsection{Arnold's scheme: Iteration}\label{sec:iteration}
Let   $d$, $\tau$, $H$, $K$, $P$, $T$, $\vae$, $\a$, $r$, $s$, $s_*$,  ${\mathsf{P}}$, $\mathsf{K}$, $\mathsf{T}$, $\torsion$, $\m$ be  as in 
Theorem~A. 
Set $K_0\coloneqq K\;,\ P_0\coloneqq P\;,\ H_0\coloneqq H$. Then, starting from $H_0$, we shall  iterate infinitely many times Lemma~\ref{lem:1bisv2}. \\
The very first step being quite different from all the others, it shall be done separately. 

\nl
Before starting, let us give some definitions\footnote{Recall the definitions of $\n$ and $\mathsf{C}_4$ given at the beginning of \S~3.1.}.
\beqano
 \m_0&\coloneqq &\m\;,\quad \torsion_0\coloneqq \torsion\;,\quad r_0\coloneqq r\;,\quad \mathsf{T}_0\coloneqq \mathsf{T}\;,\quad \mathsf{K}_0\coloneqq \mathsf{K}\;,\quad \mathsf{P}_0\coloneqq \mathsf{P}\;,\\ 
\s_0&\coloneqq& (s-s_*)/2\;,\quad \l_0\coloneqq \log\m^{-1}\;,\quad \k_0 \coloneqq 4\s_0^{-1}\l_0\;,
\\
\mathsf{C}_5   &\coloneqq& \frac{3\cdot 2^5d}{5}\;,\quad
\mathsf{C}_6 \coloneqq  \dst{\max}\left\{2^{2\n}\,,\,\mathsf{C}_5\right\}\;,\quad
\mathsf{C}_7 \coloneqq 3d\cdot 2^{6\n+2d+3}\sqrt{2}\dst\max\left\{640d^2\,,\,\mathsf{C}_4 \right\}\;,\\
\mathsf{C}_8 &\coloneqq& \left(2^{-d}\mathsf{C}_6\right)^{\su8}\;,\quad
\mathsf{C}_{9} \coloneqq 3\dst\max\left\{80d\sqrt{2}\,,\,\mathsf{C}_4\right\}\;,\\
 \l_*&\coloneqq&  \mathsf{C}_7\; \s_0^{-(4\n+2d+1)}\l_0^{2\n}\;\torsion^2\;,\quad
 \torsion_* \coloneqq 2^{2\n+2d+1}\;\mathsf{C}_6^2\;\torsion^2\;.\\
 \th_0&\coloneqq& \mathsf{C}_{9}\;\s_0^{-2(\n+d)-1}\m_0\;\torsion_0\;,\qquad  \mathsf{P}_1\coloneqq \frac{\th_0\mathsf{P}_0}{\vae\, }\;.
\eeqano
We also set, for $j\ge 0$:
\begin{align*}
& \dst\s_j\coloneqq \frac{\s_0}{2^j}\;,\quad s_{j+1}\coloneqq s_j-\s_j=s_*+\frac{\s_0}{2^j}\;,\quad \bar s_{j}\coloneqq s_j-\frac{2\s_i}{3}\;, \quad \k_j \coloneqq 4^j\k_0\;,\\ 
&  \mathsf{K}_{j+1}\coloneqq \mathsf K_0\dst\prod_{k=0}^{j}(1+\frac{\s_k}{3})\le \mathsf K_0\ex^{\frac{2\s_0}{3}}\le \mathsf{K}_0\sqrt{2}\;, \ \mathsf{T}_{j+1}\coloneqq \mathsf{T}_0\dst\prod_{k=0}^{j}(1+\frac{\s_k}{3})\le \mathsf{T}_0\sqrt{2}\,,\\
& r_{j+1}\coloneqq \su2\min\left\{\frac{\a}{2d\sqrt{2}\mathsf{K}_0\k_j^{\n}}\,,\, \frac{5}{48d}\frac{ r_j}{\torsion_0} \right\}\;,\quad \mathsf{W}_j\coloneqq \diag\left(\max\left\{\frac{\mathsf{K_j}}{{\a}}\;,\frac{1}{r_j}\right\}\;\uno_d\,,\uno_d\right)\,,\\ 
&\mathsf{L}_j\coloneqq \mathsf{P}_i\dst\max\left\{\frac{80d\sqrt{2}\;\mathsf{T}_0\;\torsion_0  }{r_j^2}\s_j^{-(\n+d)}\,,\,\mathsf{C}_4 \max\left\{1,\frac{\a}{r_j\mathsf{K}_j}\right\}\frac{ \mathsf{K}_0}{\a^2}\s_j^{-2(\n+d)}\right\}\,.
\end{align*}

\Giu
Observe that 
$$
\mathsf{W}_0=\diag\left(\mathsf{K}\a^{-1}\uno_d,\uno_d\right)\ , \quad
s_{j}\downarrow s_*\,,  \quad r_{j}\downarrow 0 \,,\quad  \ex\;\m_0\le \th_0\;.
$$
Note, also,  that, since $\th_0$ is proportional to $\e$, {\sl $\mathsf{P}_1$ is independent of $\e$.}

\subsubsection{First step}
\noi

\lem{frstStep}
Assume
\beq{condBisv2}
\a\le \frac{r_0}{\mathsf{T}_0}\qquad\mbox{and}\qquad   \th_0\le 1\;.
\eeq
Then, there exist $y_1\in D_{r_0}(y_0)$ 
 and a real--analytic symplectic transformation
\beq{phijBis0v2}
\phi_0:D_{r_1,s_{1}}(y_1)\to D_{r_{0},s_{0}}(y_{0})\ ,
\eeq
such that, for $H_1\coloneqq H_{0}\circ\phi_0$ , we have
\beq{HjBis0v2}
\left\{
\begin{aligned}
& H_1\eqqcolon K_1 + \vae^{2} P_1\;,\\
& \dpr_{y_1} K_1(y_1)=\o \;,\quad \det \dpr_{y_1}^2 K_1(y_1)\neq0  \;
\end{aligned}
\right.
\eeq\
and
\begin{align}
&|y_1-y_0| \le 
\frac{8\vae\, \mathsf{T}_0 \mathsf{P}_0}{r_0}\,,\label{y1y0diST}\\
& \|K_1\|_{r_1/4,y_1}\le \mathsf{K}_1\;,\qquad \|T_1\|\le  \mathsf{T}_1\;,\qquad T_1\coloneqq \dpr_{y_1}^2 K_1(y_1)^{-1}\;,\label{estfin2Bis0000v2}\\
& 
\vae^{2}\|P_1\|_{r_1,s_1,y_1}\le \vae^2 \mathsf{P}_1\;,\label{estfin2Bis00011v2}\\
&
\max\big\{\|\mathsf{W}_0(\phi_0-\id)\|_{r_1,s_1,y_1}\;,\, \|\dpr_x\pi_2(\phi_{0}-\id)\|_{s_1}\big\}
\le d^{-2}\s_0^{d-1}\;\vae\, {\mathsf L}_0\;.\label{estfin2Bis010v2}
\end{align}
\elem
\proof
Since
\beq{kp08}
\k_0\geby{condBisv2}4\s_0^{-1}\ge 8
\eeq
and
$$
\frac{\a}{2d\sqrt{2}\mathsf{K}_0k_0^{\n}}\overset{\equ{condBisv2}+\equ{kp08}}{\le} \frac{1}{2d\cdot 8^{\n}\sqrt{2}\mathsf{K}_0}\frac{r_0}{\mathsf{T}_0}<\frac{5}{48d}\frac{ {r}_0}{\torsion_0}\,,
$$
we get 
\beq{r1cal}
 r_1=\su2\min\left\{\frac{\a}{2d\sqrt{2}\mathsf{K}_0\k_0^{\n}}\,,\, \frac{5}{48d}\frac{ {r}_0}{\torsion_0} \right\}=\frac{\a}{4d\sqrt{2}\mathsf{K}_0\k_0^{\n}}\;.
\eeq
Thus,
\begin{align}
\vae\, \mathsf L_0 (3 \sigma_0^{-1})&\le 3\vae\, \mathsf{P}_0\dst\max\left\{\frac{80d\sqrt{2}\;\mathsf{T}_0\;\torsion_0 }{r_0^2}\s_0^{-(\n+d)}\,,\,\mathsf{C}_4 \max\left\{1,\frac{\a}{r_0\mathsf{K}_0}\right\}\frac{ \mathsf{K}_0}{\a^2}\s_0^{-2(\n+d)}\right\}\s_0^{-1}\nonumber\\
                                    &\le 3\dst\max\left\{80d\sqrt{2}\;\torsion_0\frac{ \a \;\mathsf{T}_0}{{r}_{0}}\frac{ \a }{r_0 \mathsf{K}_0}\,,\,\mathsf{C}_4 \max\left\{1,\frac{\a}{r_0\mathsf{K}_0}\right\}\right\}\s_0^{-2(\n+d)-1}\frac{ \mathsf{K}_0 \vae\, \mathsf{P}_0}{\a^2}\nonumber\\
                                    &\overset{\equ{condBisv2}}{\le}3\dst\max\left\{80d\sqrt{2}\,,\,\mathsf{C}_4\right\}\s_0^{-2(\n+d)-1}\m_0\;\torsion_0
                                    =\th_0\leby{condBisv2} 1.\label{L0ifsg3}
\end{align}
Therefore, Lemma~\ref{lem:1bisv2} implies Lemma~\ref{frstStep}.
\qed
\subsubsection{Subsequent steps, iteration and convergence}
For $j\ge 1$,  define
$$
 \m_j\coloneqq \frac{\mathsf{K}_0\,\vae^{2^j} \mathsf{P}_j}{ \a^2}\,,\qquad  \mathsf{P}_{j+1}\coloneqq 	
 			\l_*\torsion_*^{j-1} \frac{\mathsf{K}_0 { \mathsf{P}_j}^2}{\a^2}\,,\qquad \th_j     \coloneqq \l_*\;\torsion_*^j\;\m_j\,.
$$
Thus, for any $j\ge1$, one has 
\beqano
\th_{j+1}&=& \l_*\;\torsion_*^{j+1}\;\m_{j+1}=\l_*\;\torsion_*^{j+1}\frac{\mathsf{K}_0\vae^{2^{j+1}} \mathsf{P}_{j+1}}{ \a^2}=\l_*\;\torsion_*^{j+1}\frac{\mathsf{K}_0 \vae^{2^{j+1}}}{ \a^2}\;\l_*\torsion_*^{j-1} \frac{\mathsf{K}_0 { \mathsf{P}_j}^2}{\a^2}\\
		&=& \left(\l_*\;\torsion_*^{j}\;\m_j\right)^2=\th_j^2\,,
\eeqano
i.e.,
$$
\th_j=\th_1^{2^{j-1}} \;.
$$
Once the first step is completed, all the following steps do not need any other condition. Actually, 
the first condition in \equ{condBisv2} is no longer necessary and the second condition needs to be strengthen merely a little bit more. To be precise, the following holds.
\lem{lem:2Bisv2} 
Assume $\equ{HjBis0v2}\div\equ{estfin2Bis00011v2}$ 
and  
\beq{condBisv2Prt}
\mathsf{C}_8\;\torsion_0^{\su8}\;\th_1< 1\;. %* warning
\eeq
Then, one can construct a sequence of symplectic transformations 
\beq{phijBisv2}
\phi_{j-1}:D_{r_j,s_{j}}(y_j)\to D_{r_{j-1},s_{j-1}}(y_{j-1})\;,\qquad j\ge 2
\eeq
so that
\beq{HjBisv2}
H_j\coloneqq H_{j-1}\circ\phi_{j-1}=: K_j + \vae^{2^j} P_j
\eeq
converges uniformly. \\
More precisely, 
$\vae^{2^{j-1}} P_{j-1}$, 
$\phi^{j-1}\coloneqq \phi_1\circ\phi_2\circ \cdots\circ \phi_{j-1}$, 
 $K_{j-1}$, $y_{j-1}$ converge uniformly on 
$\{y_*\}\times\dst\torus^d_{s_*}$ to, respectively, $0$, $\phi^*$, $K_*$, $y_*$ and  $H_1\circ\phi^*=K_*$ with $\phi^*$ real--analytic for $x\in\dst\torus^d_{s_*}$ and $\det\dpr^2_yK_*(y_*)\neq0$. 
Finally, the following estimates hold for any $i\ge 1$:
\begin{align}
& \vae^{2^i}\|P_i\|_{r_i,s_i,y_i}\le \vae^{2^{i}} \mathsf{P}_i\ ,\label{estfin2Bis00v2}\\
&|y_{i+1}-y_i| \le \frac{8\sqrt{2}\mathsf{T}_0\vae^{2^i} \mathsf{P}_i}{r_i}\ ,\label{estfin2Bis0011v2}\\
&
|\mathsf{W}_1(\phi^*-\id)|
\le \frac{2\s_0^{d}\;\th_1}{3d^2\;\torsion_*}\qquad\quad \mbox{on}\quad \{y_*\}\times\torus^d_{s_*}\label{estfin2Bis01v2}\ .
\end{align}
\elem 

%*  Add here a remark on the difference between the first step and the other steps.
\noi
\rem{steps}
Notice that $\mathsf{P}_1$ is actually {\sl  independent of $\vae$} (and, in particular, of $\log\m^{-1}$),  while $\mathsf{P}_j$ for $j\ge 2$ does
depend on $\log \m^{-1}$ through $\l_*$. This is a crucial point, which allows, at the end, to get optimal bounds on the displacement of the persistent invariant torus from the unperturbed one.
\erem

\proof
 First of all, notice that, for any $i\ge 1$,
\begin{align*}
 r_{i+1}&= \min\left\{\frac{{\a}}{4d\sqrt{2}\mathsf{K}_0\k_i^{\n}}\,,\, \frac{5}{96d}\frac{ r_i}{\torsion_0} \right\}
	   =\min\left\{\frac{{r}_1}{4^{i\n}}\,,\, \frac{5}{96d\torsion_0}{{r}_i} \right\}\\
	   &= \min\left\{\frac{{r}_1}{4^{\n i}}\,,\, \frac{5}{96d\torsion_0}\frac{{r}_1}{ 4^{\n(i-1)}}\,,\, \left(\frac{5}{96d\torsion_0}\right)^2 {{r}_{i-1}} \right\}\\
	   &\,\ \vdots \\
	   &=\min\left\{\frac{{r}_1}{4^{\n i}}\,,\, \frac{5}{96d\torsion_0}\frac{{r}_1}{ 4^{\n(i-1)}}\,,\cdots,\, \left(\frac{5}{96d\torsion_0}\right)^i {{r}_{1}} \right\}\\
	   &=\frac{{r}_{1}}{4^{\n i}}\dst\min\left\{\left(\frac{5\cdot 4^\n}{96d\torsion_0}\right)^0\,,\,\cdots\,,\, \left(\frac{5\cdot 4^\n}{96d\torsion_0}\right)^i\right\} \\
	   &= \frac{{r}_{1}}{4^{\n i}}\dst{\min}^i\left\{\frac{5\cdot 4^\n}{96d\torsion_0}\,,\,1\right\} 
	   =  r_1 \dst{\min}^i\left\{\frac{1}{2^{2\n}}\,,\,\frac{5}{96d\torsion_0}\right\}
	   =\frac{ r_1}{\mathbf{a}_1^{i}}\,,
\end{align*}
where
\beq{eqa1maj}
\mathbf{a}_1\coloneqq \dst{\max}\left\{2^{2\n}\,,\,\frac{96d\torsion_0}{5}\right\}\le \dst{\max}\left\{2^{2\n}\,,\,\frac{96d}{5}\right\}\cdot\torsion_0=\mathsf{C}_6\;\torsion_0\;.
\eeq
\noi
For a given $j\ge 2$, let $(\mathscr{P}^j)$ be the following assertion: \\
{\sl  there exist
$j-1$ symplectic transformations\footnote{Compare \equ{phiokBis0v2}.} 
\beq{bes06v2}
\phi_{i}:D_{r_{i+1},s_{i+1}}(y_{i+1})\to D_{2r_i/3, \bar s_i}(y_i),\quad \mbox{for}\quad 1\le i\le j-1,
\eeq
 and $j-1$  Hamiltonians $H_{i+1}=H_i\circ\phi_{i}=K_{i+1}+\vae^{2^{i+1}} P_{i+1}$ real--analytic on $D_{r_{i+1},s_{i+1}}(y_{i+1})$ such that, for any $1\le i\le j-1$,
\beq{bbbBisv2}
\left\{
\begin{array}{l}
\|\dpr_y^2 K_i\|_{r_i,y_i}\le \mathsf{K}_i\,,\quad
\|T_i\|\le \mathsf{T}_i\,,\quad 
\dpr_{y} K_i(y_i)=\o \;,\quad \dpr_{y}^2 K_i(y_i)\neq0\,,\ \\  \ \\ 
\|P_{i}\|_{r_{i},s_{i},y_{i}}\le  \mathsf{P}_{i}\,,\quad
 \k_i\ge 4\s_i^{-1} \log\left(\s_i^{2\n+d}\m_i^{-1}\right)\,,\quad
\vae^{2^i} \mathsf{L}_i \le \frac{\sigma_i}{3}
\  
\end{array}\right.
\eeq
\noi
and
\beq{C.1Bisv2}
\left\{
\begin{aligned}
&\dpr_{y} K_{i+1}(y_{i+1})=\o \;,\quad \dpr_{y}^2 K_{i+1}(y_{i+1})\neq0\;,\quad
|y_{i+1}-y_i|\le \frac{8\sqrt{2}\mathsf{T}_0\vae^{2^i} \mathsf{P}_i}{r_i}\;,\\ \ \\ 
&\|T_{i+1}\|\le \|T_i\|+\mathsf T_i\vae^{2^i}\mathsf L_i\;,\quad
\|K_{i+1}\|_{r_{i+1},y_{i+1}}\le \|K_i\|_{r_i,y_{i}}+\vae^{2^i}\mathsf{P}_i \;,\\ \ \\  
&\|\dpr_y^2K_{i+1}\|_{r_{i+1},y_{i+1}}\le \|\dpr_y^2K_i\|_{r_i,y_{i}}+\mathsf K_i\vae^{2^i}\mathsf L_i \;,\\ \ \\
&\max\big\{\|\mathsf{W}_i(\phi_{i}-\id)\|_{r_{i+1},s_{i+1},y_{i+1}}\,,\,\|\dpr_x\pi_2(\phi_{i}-\id)\|_{s_{i+1}}\big\}\le d^{-2}\;\s_i^{d-1}\;\vae^{2^i}{\mathsf L}_i \;, \\ \ \\
& \|P_{i+1}\|_{r_{i+1},s_{i+1},y_{i+1}}\le  \mathsf{P}_i \mathsf L_i\;.
\end{aligned}
\right.
\eeq
}Assume $(\mathscr P^j)$, for some $j\ge 2$ and let us check $(\mathscr P^{j+1})$. Fix  $1\le i\le j-1$. Then, 
$$
\|\dpr_y^2K_{i+1}\|_{r_{i+1},y_{i+1}}\leby{C.1Bisv2} \|\dpr_y^2K_i\|_{r_i,y_{i}}+\mathsf K_i\vae^{2^i}\mathsf L_i\leby{bbbBisv2} \mathsf K_i+\mathsf K_i\frac{\s_i}{3}=\mathsf K_{i+1}<\mathsf K_0\sqrt{2}
$$
and, similarly,
$$
\|T_{i+1}\|\le \mathsf{T}_{i+1},
$$
which prove the two first relations in \equ{bbbBisv2} for $i=j$. Also
\beq{alfhtrikpi}
\frac{\a}{r_{i}\mathsf{K}_{i}}> \frac{\a}{r_1\mathsf{K}_0\sqrt{2}}=4d\k_0^{\n}\gtby{kp08}1\;,
\eeq
so that
\begin{align*}
\vae^{2^i} \mathsf L_i (3 \sigma_i^{-1})&= 3 \vae^{2^i}\mathsf{P}_i\dst\max\left\{\frac{80d\sqrt{2}\mathsf{T}_0\torsion_0  }{r_i^2}\s_i^{-(\n+d)}\,,\,\mathsf{C}_4 \max\left\{1,\frac{\a}{r_i\mathsf{K}_i}\right\}\frac{ \mathsf{K}_0}{\a^2}\s_i^{-2(\n+d)}\right\}\sigma_i^{-1}\\
                                          &\leby{alfhtrikpi} 3\vae^{2^i}\mathsf{P}_i\dst\max\left\{\frac{80d\sqrt{2}\mathsf{T}_0\torsion_0  }{r_i^2}\,,\,\mathsf{C}_4 \frac{1}{\a r_i}\right\}\s_i^{-2(\n+d)-1}\\
                                          &= 3\dst\max\left\{80d\sqrt{2}\mathsf{T}_0\torsion_0 \frac{ \a }{ {r}_{i}}\,,\,\mathsf{C}_4 \right\}\s_i^{-2(\n+d)-1}\frac{\vae^{2^i}\mathsf{P}_i}{ \a r_i}\\
                                          &= 3\dst\max\left\{640d^2\torsion_0^2\mathbf{a}^{i-1}\k_0^{\n}\,,\,\mathsf{C}_4 \right\}\s_i^{-2(\n+d)-1}\frac{\vae^{2^i}\mathsf{P}_i}{ \a^2}4d\sqrt{2}\mathsf{K}_0\k_0^{\n}\mathbf{a}^{i-1}\\
                                          &\leby{kp08} 12d\sqrt{2}\dst\max\left\{640d^2\,,\,\mathsf{C}_4 \right\}\s_i^{-2(\n+d)-1}\frac{\mathsf{K}_0\vae^{2^i}\mathsf{P}_i}{ \a^2}\torsion_0^2\mathbf{a}^{2(i-1)}\k_0^{2\n}\\
                                          &\leby{eqa1maj} 12d\sqrt{2}\dst\max\left\{640d^2\,,\,\mathsf{C}_4 \right\}\s_i^{-2(\n+d)-1}\frac{\mathsf{K}_0\vae^{2^i}\mathsf{P}_i}{\a^2}\torsion_0^{2i}\mathsf{C}_6^{2(i-1)}\k_0^{\n}\\
                                          &= 3d\cdot 2^{6\n+2d+3}\sqrt{2}\dst\max\left\{640d^2\,,\,\mathsf{C}_4 \right\}\s_0^{-(4\n+2d+1)}\left(2^{2\n+2d+1}\mathsf{C}_6^2\torsion_0^2\right)^{i-1}\cdot\\
                                          &\quad \cdot\frac{\mathsf{K}_0\vae^{2^i}\mathsf{P}_i}{ \a^2}\left(\log\m_0^{-1}\right)^{2\n}\torsion_0^2\\
                                          & 
                                          \le\mathsf{C}_7 \s_0^{-(4\n+2d+1)}\left(\log\m_0^{-1}\right)^{2\n}\torsion_0^2\:\torsion_*^{i-1}\frac{\mathsf{K}_0\vae^{2^i}\mathsf{P}_i}{\a^2}
                                          =\l_*\; \torsion_*^{i-1}\m_i
                                           = \frac{\th_i}{\torsion_*}
                                          = \frac{\th_1^{2^{i-1}}}{\torsion_*}\\
&                                          \leby{condBisv2Prt}\su{\torsion_*}<1\;.
\end{align*}
 Moreover,
 $$
 \vae^{2^{i}}\mathsf L_i<\l_*\; \torsion_*^{i-1}\m_i\;.
 $$ 
Thus, by last relation in \equ{C.1Bisv2}, for any $1\le i\le j-1$,  
 $$
 \vae^{2^{i+1}}\|P_{i+1}\|_{r_{i+1},s_{i+1},y_{i+1}}\le \vae^{2^{i}}\mathsf{L}_i\;\vae^{2^{i}} \mathsf{P}_i<\l_*\torsion_*^{i-1} \; \m_i \;\vae^{2^{i}} \mathsf{P}_i= \vae^{2^{i+1}} \mathsf{P}_{i+1}\;,
 $$ 
 which proves the fourth relation in \equ{bbbBisv2} for $i=j$. Furthermore, by exactly the same computation as above, one gets
 $$
 \vae^{2^{i+1}} \mathsf L_{i+1} (3 \sigma_{i+1}^{-1})\le \frac{\th_{i+1}}{\torsion_*}=\frac{\th_1^{2^{i}}}{\torsion_*}<1\ ,
 $$
 which proves the last relation in \equ{bbbBisv2} for $i=j$.  It remains only to check that
the fifth relation in \equ{bbbBisv2} holds as well for $i=j$ in order to apply Lemma~\ref{lem:1bisv2} to $H_i$, $1\le i\le j$ and get \equ{C.1Bisv2} and, consequently, $(\mathscr P^{j+1})$. In fact,  we have\footnote{Notice that $\left(\log t\right)^{2s}\le t^{1/2}\;,\quad\forall\;t\ge \ex,\quad\forall\;s\ge 1/4,$ so that $\m_0(\log\m_0^{-1})^{2\n}\leby{condBisv2Prt} \sqrt{\m_0}\le \ex^{-1/2}<1$, which in turn proves the {\it r.h.s.} inequality in \equ{th1th0m0}.}
 \beq{th1th0m0}
\l_*\;\torsion_*\;\m_0^2<\l_*\;\torsion_*\;\m_0\;\th_0=\th_1\le \mathsf{C}_7\;\s_0^{-(4\n+2d+1)}\;\torsion_*\;\torsion_0^2\;\th_0\;,
 \eeq
 so that
 \begin{align*}
 4\s_j^{-1} \log\left(\s_j^{2\n+d}\m_j^{-1}\right)&\le 4\s_j^{-1} \log\left(\m_j^{-1}\right)
         = 4\s_j^{-1} \log\left({\l_*\torsion_*^{j}}\th_1^{-2^{j-1}}\right)\\
         &\leby{th1th0m0} 4\s_j^{-1} \log\left({\l_*\torsion_*^{j}} ({\l_*\torsion_*}\m_0^2)^{-2^{j-1}}\right)
         \le 4\s_j^{-1}\log\left(\m_0^{-2^{j}}\right)\\
         &= 4^j\cdot 4\s_0^{-1} \log\left(\m_0^{-1}\right)
         = \k_j \;.
 \end{align*}
\noi
To finish the proof of the induction, i.e., to  construct an {\sl infinite sequence} of Arnold's transformations satisfying \equ{bbbBisv2} and  \equ{C.1Bisv2} {\sl for all $i\ge 1$}, one needs only to check $(\mathscr P^{2})$. Thanks to\footnote{Observe that for $j=2$, $i=1$.} $\equ{HjBis0v2}\div\equ{estfin2Bis00011v2}$, we just need to check the two last inequalities in $\equ{bbbBisv2}_{i=1}$. But, in fact, this is contained in the above computation. Then, we apply Lemma~\ref{lem:1bisv2} to $H_1$ to get $\equ{bes06v2}_{i=1}$ and  $\equ{C.1Bisv2}_{i=1}$, which achieves the proof of $(\mathscr P^{2})$.\\
\nl
Next, we prove that $\phi^j$ is convergent by proving that it is a Cauchy sequence. For any $j\ge 4$, we have, using again Cauchy's estimate (and noting that $2^{i-1}\ge i,\,\forall\; i\ge 0$),
\beqano
\|\mathsf{W}_{j-1}(\phi^{j-1}-\phi^{j-2})\|_{r_j,s_j,y_j}&=&\|\mathsf{W}_{j-1}\phi^{j-2}\circ\phi_{j-1}-\mathsf{W}_{j-1}\phi^{j-2}\|_{r_j,s_j,y_j}\\
           &\leby{bes06v2}& \|\mathsf{W}_{j-1}D\phi^{j-2}\mathsf{W}_{j-1}^{-1}\|_{2r_{j-1}/3, s_{j-1},y_{j-1}}\, \|\mathsf{W}_{j-1}(\phi_{j-1}-\id)\|_{r_j,s_j,y_j}\\
           &\leby{C.1Bisv2}&  \max\left(r_{j-1}\frac{3}{r_{j-1}},\frac{3}{2\s_{j-1}}\right)    \|\mathsf{W}_{j-1}\phi^{j-2}\|_{r_{j-1}, s_{j-1},y_{j-1}} \times\\
           &&\qquad \times \|\mathsf{W}_{j-1}(\phi_{j-1}-\id)\|_{r_j,s_j,y_j}\\
           &=&  \frac{3}{2\s_{j-1}}   \|\mathsf{W}_{j-1}\phi^{j-2}\|_{r_{j-1}, s_{j-1},y_{j-1}} \, \|\mathsf{W}_{j-1}(\phi_{j-1}-\id)\|_{r_j,s_j,y_j}\\
           &\le & \frac{1}{2}    \|\mathsf{W}_{j-1}\phi^{j-2}\|_{r_{j-1}, s_{j-1},y_{j-1}} \cdot \s_{j-1}^d\left(\vae^{2^{j-1}}{\mathsf{L}}_{j-1}3\s_{i-1}^{-1}\right)\\
           &\le & \frac{1}{2}    \|\mathsf{W}_{j-1}\phi_1\|_{r_{2}, s_2,y_{2}} \cdot \s_{j-1}^d\;\th_{j-1}\\
           &\le &  \frac{1}{2}\left(\dst\prod_{i=1}^{j-2}\|\mathsf{W}_{i+1}\mathsf{W}_{i}^{-1}\| \right)\|\mathsf{W}_{1}\phi_1\|_{r_{2}, s_2,y_{2}} \cdot \s_{j-1}^d\;\th_{j-1}\\
           &\eqby{alfhtrikpi}& \frac{1}{2}\left(\dst\prod_{i=1}^{j-2}\frac{r_i}{r_{i+1}} \right)\|\mathsf{W}_{1}\phi_1\|_{r_{2}, s_2,y_{2}} \cdot \s_{j-1}^d\;\th_{j-1}\\
           &=& \frac{r_1}{2r_{j-1}}\|\mathsf{W}_{1}\phi_1\|_{r_{2}, s_2,y_{2}} \cdot \s_{j-1}^d\;\th_{j-1}\\
           &\leby{eqa1maj}& \su2\s_{3}^d\;(\mathsf{C}_6\;\torsion_0)^2\;\|\mathsf{W}_{1}\phi_1\|_{r_{2}, s_2,y_{2}} \cdot \left(2^{-d}\mathsf{C}_6\;\torsion_0\right)^{j-4}\cdot \;\th_1^{2^{j-2}}\\
           &\le& \su2\s_{3}^d\;(\mathsf{C}_6\;\torsion_0)^2\;\|\mathsf{W}_{1}\phi_1\|_{r_{2}, s_2,y_{2}} \cdot \left(2^{-d}\mathsf{C}_6\;\torsion_0\right)^{2^{j-5}}\cdot \;\th_1^{2^{j-2}}\\
           &=& \su2\s_{3}^d\;(\mathsf{C}_6\;\torsion_0)^2\;\|\mathsf{W}_{1}\phi_1\|_{r_{2}, s_2,y_{2}} \cdot \left(\left(2^{-d}\mathsf{C}_6\;\torsion_0\right)^{\su8} \th_{1}\right)^{2^{j-2}}\\
           &=& \su2\s_{3}^d\;(\mathsf{C}_6\;\torsion_0)^2\;\|\mathsf{W}_{1}\phi_1\|_{r_{2}, s_2,y_{2}} \cdot \left(\mathsf{C}_8\;\torsion_0^{\su8}\; \th_{1}\right)^{2^{j-2}}  \;.
\eeqano
\noi
Therefore, for any $n\ge 1,\, j\geq 0$,
\begin{align*}
\|\mathsf{W}_{1}(\phi^{n+j+1}-\phi^n)\|_{r_{n+j+2},s_{n+j+2},y_{n+j+2}}&\leq  \sum_{i=n}^{n+j}\|\mathsf{W}_{1}(\phi^{i+1}-\phi^i)\|_{r_{i+2},s_{i+2},y_{i+2}}\\
&\le \sum_{i=n}^{n+j}\left(\dst\prod_{k=1}^{i}\|\mathsf{W}_{k}\mathsf{W}_{k+1}^{-1}\| \right)\|\mathsf{W}_{i+1}(\phi^{i+1}-\phi^i)\|_{r_{i+2},s_{i+2},y_{i+2}}\\
&\eqby{alfhtrikpi} \sum_{i=n}^{n+j}\dst\prod_{k=1}^{i}\max\left\{1,\frac{r_{k+1}}{r_k} \right\}\|\mathsf{W}_{i+1}(\phi^{i+1}-\phi^i)\|_{r_{i+2},s_{i+2},y_{i+2}}\\
&= \sum_{i=n}^{n+j}\|\mathsf{W}_{i+1}(\phi^{i+1}-\phi^i)\|_{r_{i+2},s_{i+2},y_{i+2}}\\
&\le \su2\s_{3}^d\;(\mathsf{C}_6\;\torsion_0)^2\;\|\mathsf{W}_{1}\phi_1\|_{r_{2}, s_2,y_{2}}\cdot\sqrt{\vae\,} \dst\sum_{i=n}^{n+j} \left(\mathsf{C}_8\;\torsion_0^{\su8}\; \th_{1}\right)^{2^{i}}\;.
\end{align*}
Hence, by \equ{condBisv2Prt}, $\phi^j$ converges uniformly on $\{y_*\}\times\torus^d_{s_*}$ to some $\phi^*$, which is then real--analytic map in $x\in\torus^d_{s_*}$.

\nl
To estimate $|\mathsf{W}_0(\phi^*-\id)|$ on $\{y_*\}\times\torus^d_{s_*}$, observe that
, for $i\ge 1$,
$$\s_{i}^d\;\vae^{2^i}\mathsf L_i\le \frac{\s_0^{d+1}}{3 \cdot 2^{i(d+1)}}\ \frac{\th_1^{2^{i-1}}}{\torsion_*} \le \frac{\s_0^{d+1}}{3 \cdot 2^{(d+1)i} \torsion_*}\th_1^{{i}}= \frac{(2\s_0)^{d+1}}{3\torsion_*} \Big(\frac{\th_1}{2^{d+1}}\Big)^{i+1}$$
and therefore 
$$\dst\sum_{i\ge 1}  \s_{i}^d\;\vae^{2^i}\mathsf L_i\le \frac{(2\s_0)^{d+1}}{3\torsion_*}\sum_{i\ge 1}\Big(\frac{\th_1}{2^{d+1}}\Big)^{i}\le \frac{2\s_0^{d+1}\;\th_1}{3\;\torsion_*}
\ .$$ 
Moreover, for any $i\ge 1$,
\begin{align*}
\|\mathsf{W}_1(\phi^i-\id)\|_{r_{i+1},s_{i+1},y_{i+1}}&\le \|\mathsf{W}_1(\phi^{i-1}\circ\phi_i-\phi_i)\|_{r_{i+1},s_{i+1},y_{i+1}}+\|\mathsf{W}_1(\phi_i-\id)\|_{r_{i+1},s_{i+1},y_{i+1}}\\
&\le \|\mathsf{W}_1(\phi^{i-1}-\id)\|_{r_{i},s_{i},y_{i}}+ (\dst\prod_{j=0}^{i-1}\|\mathsf{W}_{j}\mathsf{W}_{j+1}^{-1}\|) \|\mathsf{W}_{i}(\phi_i-\id)\|_{r_{i+1},s_{i+1},y_{i+1}}\\
&= \|\mathsf{W}_1(\phi^{i-1}-\id)\|_{r_{i},s_{i},y_{i}}+ \|\mathsf{W}_{i}(\phi_i-\id)\|_{r_{i+1},s_{i+1},y_{i+1}}\\
&= \|\mathsf{W}_1(\phi^{i-1}-\id)\|_{r_{i},s_{i},y_{i}}+ \|\mathsf{W}_{i}(\phi_i-\id)\|_{r_{i+1},s_{i+1},y_{i+1}}\\
&\le \|\mathsf{W}_1(\phi^{i-1}-\id)\|_{r_{i},s_{i},y_{i}}+d^{-2}\s_{i}^{d-1}\;\vae^{2^{i}}{\mathsf{L}}_{i}\ ,
\end{align*}
which iterated yields
\begin{align*}
\|\mathsf{W}_1(\phi^i-\id)\|_{r_{i+1},s_{i+1},y_{i+1}}&\le d^{-2}\dst\sum_{k\ge 1}\s_{k}^{d-1}\; \vae^{2^k}{\mathsf{L}}_k\le \frac{2\s_0^{d}\;\th_1}{3d^2\;\torsion_*}
\,.
\end{align*}
Therefore, taking the limit over $i$ completes the proof of \equ{estfin2Bis01v2} and hence of Lemma~\ref{lem:2Bisv2}.
\qed
\subsection{Conclusion}
We can now complete the proof of Theorem~A. Let
\begin{align*}
&\mathsf{C}_{10} \coloneqq \left(2^{-(4\n+2d+1)}  +2\mathsf{C}_7\right)\mathsf{C}_9/(3d^2)\;,\quad
\mathsf{C}_{11} \coloneqq \frac{1}{2^{5\n+3d-2}}+\frac{\mathsf{C}_7\;\mathsf{C}_9}{3\cdot 5\cdot 2^{\n+2}\cdot d^2\cdot \sqrt{2}}\;,\\
&\mathsf{C}_{12} \coloneqq 2^{2\n+2d+1}\;\mathsf{C}_6^2\;\mathsf{C}_7\;\mathsf{C}_8\;\mathsf{C}_9\;,\quad
\mathsf{C}_{13} \coloneqq \mathsf{C}_{10}+2^{-(\n+1)}\;\mathsf{C}_{11}\;\quad, \mathsf{C}_{14}\coloneqq 2^{2(3\n+2d+1)} \mathsf{C}_{12}\;,\\
&\mathsf{C}_{15}\coloneqq 18d^3+70\;,\quad\mathsf{C} \coloneqq 2^{6\t+3d+8} \mathsf{C}_{13}\;,\quad 
\mathsf{C}_*\coloneqq  \max\left\{(4\n\ex^{-1})^{8\n/3}\mathsf{C}_{14}^{2/3},\;  2\mathsf{C}_{15}\mathsf{C}\right\}\;.
\end{align*}
Observe that
\beq{stIneq00}
(\log t)^{4\n}\le (4\n\ex^{-1})^{4\n}\sqrt{t}\;,\qquad \forall\; t>1.
\eeq
Then,
\begin{align*}
\mathsf{C}_8\;\torsion_0^{\su8}\;\th_1&=\mathsf{C}_{14}\;\torsion^{41/8}\;(s-s_*)^{-2(3\n+2d+1)}\;\m^2(\log \m^{-1})^{2\n}\\
  &\leby{stIneq00} (4\n\ex^{-1})^{4\n}\;\mathsf{C}_{14}\;\torsion^{41/8}\;\m^{3/2}\;(s-s_*)^{-2(3\n+2d+1)}\\
  &< \left(\mathsf{C}_*\;\torsion^4\;(s-s_*)^{-(6\n+3d+2)} \;\m\right)^{3/2}\\
  &\leby{smcondwhL} 1
\end{align*}
and
$$
\th_0<\mathsf{C}_*\;\torsion^4\;(s-s_*)^{-(6\n+3d+2)} \;\m\leby{smcondwhL} 1.
$$
Hence, \equ{smcondwhL} implies the smallness conditions \equ{condBisv2} and \equ{condBisv2Prt}. Therefore, Lemma~\ref{frstStep} and \ref{lem:2Bisv2} hold. Now, set $\phi_*\coloneqq \phi_0\circ \phi^*$ and observe that, uniformly on $\{y_*\}\times \torus^d_{s_*}$,
\begin{align*}
|\mathsf{W}_0(\phi_*-\id)|&\le |\mathsf{W}_0(\phi_0\circ \phi^*-\phi^*)|+|\mathsf{W}_0(\phi^*-\id)|\\
&\le \|\mathsf{W}_0(\phi_0-\id)\|_{r_1,s_1,y_1}+\|\mathsf{W}_0\mathsf{W}_1^{-1}\|\;|\mathsf{W}_1(\phi^*-\id)|\\
&\le \su{d^2}\s_{0}^d\;\vae\,{\mathsf{L}}_0+\frac{2\s_0^{d}}{3d^2\torsion_*}\;\th_1
\stackrel{\equ{L0ifsg3}+\equ{th1th0m0}}{\le} \frac{\s_0^{d}}{3 d^2}\;\th_0 +\frac{2\s_0^{d}}{3d^2\torsion_*}\mathsf{C}_7\;\s_0^{-(4\n+2d+1)}\;\torsion_*\;\torsion_0^2\;\th_0\\
&\le \left(\frac{1}{3d^2 {2^{4\n+2d+1}} } +\frac{2\mathsf{C}_7}{3d^2}\right)\s_0^{-(4\n+d+1)}\;\torsion_0^2\;\th_0
=\mathsf{C}_{10} \;\s_0^{-(6\n+3d+2)}\;\torsion_0^3\;\m_0\eqqcolon \g
\;.
\end{align*}
 Moreover, for any $i\ge 1$,
\begin{align*}
|y_{i}-y_0|&\le \dst\sum_{j=0}^{i-1}|y_{j-1}-y_j|
 	\stackrel{\equ{y1y0diST}+\equ{estfin2Bis0011v2}}{\le}\frac{8\mathsf{T}_0\vae\, \mathsf{P}_0}{r_0}+\dst\sum_{j=1}^{i-1}\frac{8\sqrt{2}\mathsf{T}_0\vae^{2^j} \mathsf{P}_j}{r_j}\\
 	&\le \frac{8\mathsf{T}_0\vae\, \mathsf{P}_0}{r_0}+\dst\sum_{j\ge 1}\frac{r_j}{10d\torsion_0}\s_j^{\n+d}\vae^{2^j}\mathsf{L}_j
 	\le \frac{8\mathsf{T}_0\vae\, \mathsf{P}_0}{r_0}+\frac{r_1}{10d\torsion_0}\s_1^{\n}\dst\sum_{j\ge 1}\s_j^{d}\vae^{2^j}\mathsf{L}_j\\
 	&\le \frac{8\mathsf{T}_0\vae\, \mathsf{P}_0}{r_0}+\frac{r_1}{10d\torsion_0}\s_1^{\n}\;\frac{2\s_0^{d+1}\th_1}{3\torsion_*}
 	\stackrel{\equ{smcondwhL}+\equ{th1th0m0}}{\le}\mathsf{C}_{11}\;\s_0^{-(5\n+3d+1)}\;\torsion_0^2\;\frac{\vae\,\;\mathsf{P}_0}{\a}\;,
\end{align*}
and then passing to the limit, we get
$$
|y_{*}-y_0|\le  \mathsf{C}_{11}\;\s_0^{-(5\n+3d+1)}\;\torsion_0^2\;\frac{\vae\,\;\mathsf{P}_0}{\a}\;.
$$
Thus, the triangle inequality gives
$$
\dst\sup_{ \torus^d_{s_*}}|\mathsf{W}_0(\phi_*-\phi_{\rm e})|\le \mathsf{C}_{13}\;\s_0^{-(6\n+3d+2)}\;\torsion_0^3\;\m_0\;,
$$
which proves the bounds on $\|u_*\|$ and $\|v_*\|$ in \equ{est}. 
 Let us now prove the bound on $\dpr_x u_*$   in \equ{est}.
 Set\\
$$\tilde{u}_j\coloneqq \dpr_x\pi_2(\phi_j-\id),\qquad U^j\coloneqq\dpr_x\pi_2\phi_0\circ\phi^j=(\uno_d+\tilde{u}_0)\circ\cdots\circ(\uno_d+\tilde{u}_j).
$$
 Then, for any $j\ge 0$, we have
 $$
 \|U^j\|_{s_{i+1}}\le (1+\|\tilde{u}_0\|_{s_0})\cdots(1+\|\tilde{u}_j\|_{s_j})\stackrel{\equ{estfin2Bis010v2}+\equ{C.1Bisv2}}{\le}\exp\left(d^{-2}\dst\sum_{k\ge 0}\s_{k}^{d-1}\; \vae^{2^k}{\mathsf{L}}_k\right)\le \ex^\g\;,
 $$
 so that
 $$
 \|U^{j+1}-U^j\|_{s_*}=\|U^{j}(\uno_d+\tilde{u}_{j+1})-U^j\|_{s_*}\le \|U^j\|_{s_{j+1}}\|\tilde{u}_{j+1}\|_{s_{j+1}} \stackrel{\equ{estfin2Bis010v2}+\equ{C.1Bisv2}}{\le}\ex^\g d^{-2}\dst\s_{j+1}^{d-1}\; \vae^{2^{j+1}}{\mathsf{L}}_{j+1},
 $$
 which implies
 $$
  \|U^{j}-\uno_d\|_{s_*}\le \ex^\g d^{-2}\sum_{k\ge 0}\s_{k}^{d-1}\; \vae^{2^{k}}{\mathsf{L}}_{k}\le \g\ex^\g\leby{smcondwhL} \ex\;\g\leby{smcondwhL}\su2
 $$
 and then letting $j\rightarrow\infty$, we get the estimate on $\dpr_x u_*$.
  \qed
  
  \rem{rem:final}
As it is easy to check,  Theorem~A  holds under the milder condition $\m\le \m_\sharp$ where
\beqano
\m_\sharp 		  &\coloneqq\max\big\{\m>0\;: &\mathsf{C}_{14}\;\torsion^{\frac{41}{8}}\;(s-s_*)^{-2(3\n+2d+1)}\;\m^2\left(\log\m^{-1} \right)^{2\n}\le 1\  ,\\
               &   \ \qquad {\rm and} & 2 \mathsf{C} \;(s-s_*)^{-(6\n+3d+2)}\;\torsion^3\;\m\;\exp\left(\mathsf{C} \;(s-s_*)^{-(6\n+3d+2)}\;\torsion^3\;\m\right)\le 1	\big\}\;.
\eeqano
Notice that $\m_*<\m_\sharp$.\\
Indeed,  condition 
$$
\mathsf{C}_{14}\;\torsion^{\frac{41}{8}}\;(s-s_*)^{-2(3\n+2d+1)}\;\m^2\left(\log\m^{-1} \right)^{2\n}\le 1
$$
guaranties the convergence of Arnold's scheme, while  condition
$$
2 \mathsf{C} \;(s-s_*)^{-(6\n+3d+2)}\;\torsion^3\;\m\;\exp\left(\mathsf{C} \;(s-s_*)^{-(6\n+3d+2)}\;\torsion^3\;\m\right)\le 1,
$$
ensures that the Torus $\mathcal{T}_{\o,\vae}$ is a Lagrangian  graph (over the ``angle'' variables).

\erem

\Giu
{\bf Acknowledgements} L.C. has been supported by the ERC grant HamPDEs under FP7 n.
306414 and the PRIN national grant ``Variational Methods in Analysis, Geometry and Physics''.
The authors are indebted with an anonymous referee for valuable suggestions and corrections.

\appendix
\section*{Appendix}
\addcontentsline{toc}{section}{Appendices}
\setcounter{section}{0}
\renewcommand{\thesection}{\Alph{section}}

\appA{Constants}\label{AssumpExtArnolv2}
For convenience, we collect here  the list of  constants  appearing in the proof of Theorem~A. 
\\
Recall that 
$\t\ge d-1\ge 1$ and notice that all  $\mathsf{C}_i$'s are  greater than $1$ and depend only upon $d$ and $\t$.
\beqano
\n			 &\coloneqq& \t+1\;,\\
\mathsf{C}_0 &\coloneqq& 4\sqrt{2}\left(\frac{3}{2}\right)^{2\n+d}\dst\int_{\rn} \left( |y|_1^{\n}+|y|_1^{2\n}\right)\ex^{-|y|_1}dy\;,\quad
\mathsf{C}_1 \coloneqq 2\left(\frac{3}{2}\right)^{\n+d}\dst\int_{\rn} |y|_1^{\n}\ex^{-|y|_1}dy\;,\\
\mathsf{C}_2 &\coloneqq& 2^{3d}d\;,\quad
\mathsf{C}_3 \coloneqq	\left(d^2\mathsf{C}_1^2+6d\mathsf{C}_1 +\mathsf{C}_2\right)\sqrt{2}\;,\quad
\mathsf{C}_4 \coloneqq \max\left\{6d^2\mathsf{C}_0,\,\mathsf{C}_3\right\}\;,\\
\mathsf{C}_5   &\coloneqq& \frac{3\cdot 2^5d}{5}\;,\quad
\mathsf{C}_6 \coloneqq \dst{\max}\left\{2^{2\n}\,,\,\mathsf{C}_5\right\}\;,\quad
\mathsf{C}_7 \coloneqq 3d\cdot 2^{6\n+2d+3}\sqrt{2}\dst\max\left\{640d^2\,,\,\mathsf{C}_4 \right\}\;,\\
\mathsf{C}_8 &\coloneqq& \left(2^{-d}\mathsf{C}_6\right)^{\su8}\;,\quad
\mathsf{C}_{9} \coloneqq 3\dst\max\left\{80d\sqrt{2}\,,\,\mathsf{C}_4\right\}\;,\quad
\mathsf{C}_{10} \coloneqq \left(2^{-(4\n+2d+1)}  +2\mathsf{C}_7\right)\mathsf{C}_9/(3d^2)\;,\\
\mathsf{C}_{11} &\coloneqq& \frac{1}{2^{5\n+3d-2}}+\frac{\mathsf{C}_7\;\mathsf{C}_9}{3\cdot 5\cdot 2^{\n+2}\cdot d^2\cdot \sqrt{2}}\;,\quad
\mathsf{C}_{12} \coloneqq 2^{2\n+2d+1}\;\mathsf{C}_6^2\;\mathsf{C}_7\;\mathsf{C}_8\;\mathsf{C}_9\;,\\
\mathsf{C}_{13} &\coloneqq& \mathsf{C}_{10}+2^{-(\n+1)}\;\mathsf{C}_{11}\;,\quad
\mathsf{C}_{14}\coloneqq 2^{2(3\n+2d+1)} \mathsf{C}_{12}\;,\quad
\mathsf{C}_{15}\coloneqq 18d^3+70\;,\\
\mathsf{C}\ &\coloneqq& 2^{6\t+3d+8} \mathsf{C}_{13}\;,\quad 
\mathsf{C}_*\coloneqq  \max\left\{(4\n\ex^{-1})^{8\n/3}\mathsf{C}_{14}^{2/3},\;  2\mathsf{C}_{15}\mathsf{C}\right\}\;.
\eeqano

\appB{Kolmogorov's non--degeneracy}\label{salamon}
Let
$$
\hat \m \coloneqq 2 \ex \ \mathsf{C}\;(s-s_*)^{-(6\t+3d+8)}\;\torsion^3\;\cdot \frac{\mathsf{K}\;\vae \mathsf{P}  }{\a^2}\;.
$$
Since $\|\dpr_x u_*\|_{s_*} \stackrel{\equ{est}}{\le}1/2$, then $\id+u_*$ is a diffeomorphism of $\tn$.
Letting
$$(\dpr_x(\id+u_*)(x))^{-1} \eqqcolon \uno_d+A(x)\ ,
$$
we have 
\beq{equsttiInv}
\|A\|_{s_*}\le 2\|\dpr_x u_*\|_{s_*}\leby{est} 2\hat \m\leby{smcondwhL} 1\;;\quad
\|v_*\|_{s_*}\leby{est} \frac{\a}{\mathsf{K}}\frac{\hat \m}{2\ex}\leby{smcondwhL} \frac{r}{\torsion}\frac{1}{4\ex\mathsf{C}_{15}}< \frac{r}{8}.
\eeq
In \cite{Salamon} it is proven that the map
$$
\phi(y,x)\coloneqq (y_0+v_*(x)+y+ A^{T}y,x+u_*(x)). %* warning
$$
is symplectic. Then, 
$$
H\circ \phi (y,x)= E+\o \cdot y + Q(y,x) 
$$
with:
\begin{align*}
&E= K(y_0),\quad \langle Q_{yy}(0,\cdot)\rangle=K_{yy}(y_0)+\average{M}\;,\\
& M\coloneqq \dpr^2_y\bigg(K(y_0+v_*+y+A^{T}y)-\su2 y^TK_{yy}(y_0)y\bigg)\Big|_{y=0}+\dpr^2_y(\vae P\circ\phi)\Big|_{y=0}\;,\\
& \sup_{\torus^d_{s_*}}\|K_{yy}(y_0)^{-1}M\|\stackrel{\equ{equsttiInv}}{\le} (18d^3+70)\hat \m\torsion\leby{smcondwhL} 1/2,
\end{align*}
which show that $\langle Q_{yy}(0,\cdot)\rangle$ is invertible.

\appC{Reminders}

\subsection{Classical estimates (Cauchy, Fourier)} 
\lem{Cau}  {\rm \cite{CC95}} 
Let $p\in \natural,\,r,s>0, y_0\in \cn$ and $f$ a real--analytic function $D_{r,s}(y_0)$ with 
$$
\|f\|_{r,s}\coloneqq \dst\sup_{D_{r,s}(y_0)}|f|.
$$
Then,\\
{\bf (i)} For any multi--index $(l,k)\in \natural^d\times\natural^d$ with $|l|_1+|k|_1\le p$ and for any $0<r'<r,\, 0<s'<s$,\footnote{As usual, $\dpr_y^l\coloneqq \frac{\dpr^{|l|_1}}{\dpr y_1^{l_1}\cdots\dpr y_d^{l_d}},\, \forall\, y\in\rn,\, l\in\zn $.\label{notDevPart}}
\[\|\partial_{y}^l \partial_{x}^k f\|_{r',s'}\leq p!\; \|f\|_{r,s}(r-r')^{|l|_1}(s-s')^{|k|_1}.\]
{\bf (ii)}  For any $ k\in \zn$ and any $y\in D_r(y_0)$ 
$$|f_k(y)|\leq \ex^{-|k|_1 s}\|f\|_{r,s}.
$$
\elem
\subsection{Implicit function theorem}
\lem{IFTLem} {\rm \cite{CL12}}
Let $ r,s>0,\, n,m\in \natural,\, (y_0,x_0)\in \complex^n\times\complex^m$ and\footnote{Here, $D^n_r(z_0)$ denotes the ball in $\complex^n$ centered at $z_0$ and with radius $r$.} 
\[F\colon (y,x)\in D^n_r(y_0)\times D^m_s(x_0)\subset \complex^{n+m}\mapsto F(y,x)\in\complex^n\]
be continuous with continuous Jacobian matrix $F_y$. Assume that $F_y(y_0,x_0)$ is invertible with inverse $T\coloneqq F_y(y_0,x_0)^{-1}$ 
 such that
\beq{HypIFT}
\sup_{D^n_r(y_0)\times D^m_s(x_0)}\|\uno_n-TF_y(y,x)\|\leq c<1 \quad \mbox{and}\quad \sup_{ D^m_s(x_0)}|F(y_0,\cdot)|\leq \frac{(1-c) r}{\|T\|}.
\eeq
Then, there exists a unique continuous function $ g\colon D^m_s(x_0)\to D^n_r(y_0)$ such that the following are equivalent
\begin{itemize}
\item[$(i)$] $(y,x)\in D^n_r(y_0)\times D^m_s(x_0)$ and $F(y,x)=0$;
\item[$(ii)$] $x\in D^m_s(x_0)$ and $y=g(x)$.
\end{itemize}
Moreover, $g$ satisfies
\beq{EstIFT}
\sup_{D^m_s(x_0)}|g-y_0|\leq \frac{\|T\|}{1-c}\sup_{D^m_s(x_0)}|F(y_0,\cdot)|.
\eeq
\elem
%%%%%%%%%%%%

\end{document}